\documentclass[12pt]{amsart}
\usepackage{amsthm}
\usepackage{epsf}
\usepackage{amssymb}
\usepackage{latexsym}
\usepackage{amsmath}
\newtheorem{theorem}{Theorem} 
\newtheorem{definition}[theorem]{Definition}

\newtheorem{corollary}[theorem]{Corollary}
\newtheorem{lemma}[theorem]{Lemma}
\newtheorem{example}[theorem]{Example}
\newtheorem{proposition}[theorem]{Proposition}
\newtheorem{remark}[theorem]{Remark}

\newtheorem{notation}[theorem]{Notation}

\def\R{{\mathbb R}}

\def\C{{\mathbb C}}

\def\Z{{\mathbb Z}}
\def\P{{\mathbb P}}

\def \ta{\tau}
\def \ta1{\tau_1}

\def \Dl{\Delta}
\def \dl{\delta}
\def \g{\gamma}

\def \vp{\varphi}

\def \tC{\tilde{C}}

\def \prodl{\prod\limits}

\def\bZ{\bar{Z}}

\newcommand\set[1]{{\{{#1}\}}}

\def\st{{such that }}

\newcommand\begintable[1][] {{}}

\long\def\forget#1\forgotten{}

\newif\ifXY 
\XYtrue     

\begin{document}

\renewcommand{\subjclassname}{%
       \textup{2000} Mathematics Subject Classification}

\title{Braid monodromy factorization for a non-prime $K3$ surface branch curve}

\author[M.~Amram, C.~Ciliberto, R.~Miranda, M.~Teicher]
{Amram Meirav$^1$, Ciliberto Ciro, Miranda Rick and Teicher Mina}
\stepcounter{footnote}
\footnotetext{Partially supported by the DAAD fellowship (Germany),
the Golda Meir postdoctoral fellowship (the Einstein
mathematics institute, Hebrew university, Jerusalem),
the Emmy Noether Research Institute for Mathematics
(center of the Minerva Foundation of Germany), the Excellency Center "Group
Theoretic Methods in the Study of Algebraic Varieties"
of the Israel Science Foundation, and EAGER (EU network, HPRN-CT-2009-00099).}
\address{Meirav Amram, Mathematisches Institut, Bismarck Strasse 1 1/2, Erlangen,
Germany; Einstein Institute for Mathematics, Hebrew University, Jerusalem}
\email{meirav@macs.biu.ac.il/ameirav@math.huji.ac.il}
\address{Ciro Ciliberto, Dipartimento di Matematica, Universita di Roma II,
Tor Vergata,
00123 Roma, Italy}
\email{cilibert@mat.uniroma2.it}
\address{Rick Miranda, Department of Mathematics, Colorado State University,
Fort Collins, CO 80523 USA}
\email{rick.miranda@colostate.edu}
\address{Mina Teicher, Department of Mathematics, Bar-Ilan university, 52900
Ramat-Gan, Israel}
\email{teicher@macs.biu.ac.il}

\begin{abstract}
In this paper we consider a non-prime $K3$ surface of degree $16$,
and study a specific degeneration of it,
known as the $(2,2)$-pillow degeneration, \cite{11}.

We study also the braid monodromy factorization
of the branch curve of the surface
with respect to a generic projection onto $\C\P^2$.

In \cite{pil} we compute the fundamental groups
of the complement of the branch curve
and of the corresponding Galois cover of the surface.
\end{abstract}


\date{\today}

\maketitle
\section{Overview}\label{overview}

Given a projective surface and a generic projection to the plane,
the fundamental group of the complement of the branch curve
is one of its most important invariants.
Our goal is to compute this group
and the fundamental group of the Galois cover
(which is known to be a certain quotient of the fundamental group
of the complement of the branch curve, see \cite{15}).
This goal is achieved in \cite{pil}.

In this paper we deal with a non-prime $K3$ surface of degree 16
which is embedded in $\P^9$.
In order to compute the above groups
we degenerate the surface into a union of 16 planes.
We then project it onto $\C\P^2$ to get a degenerated branch curve $S_0$,
from which one can compute the braid monodromy factorization
and the branch curve $S$ of the projection of the original $K3$ surface.
With this information we will be able
to apply the van Kampen Theorem (see \cite{20})
and the regeneration rules (see \cite{19})
to get presentations for the relevant fundamental groups.

The idea of using degenerations for these purposes appears already in
\cite{3}, \cite{10}, \cite{12}, \cite{15} and \cite{18}.
Degenerations of $K3$ surfaces, of which the one we use here is an example,
were constructed in \cite{11} and are called {\em pillow degenerations}.

The study of braids was begun by Artin in \cite{1} and \cite{2};
see also Chisini \cite{9}.
The braid monodromy technique was first presented by Moishezon-Teicher
in \cite{13}, \cite{16} and \cite{17}.
Many examples of computations of braid monodromy have been executed,
see for example in \cite{3}, \cite{5}, \cite{17}.
Also in \cite{dennis} one can find a description
of computations of braid monodromy
and the fundamental group $\pi_1(\C\P^2-S)$
of the complement of the branch curve of a surface.
In this work we encounter new type of singularities ($3$-points, see Section \ref{3pt})
which have not been handled before,
whose analysis is necessary to give the precise computations of the braid monodromy.
Moreover, since the monodromies related to $6$-points are quite hard to follow,
we present them first in a precise way algebraically,
followed then by figures which illustrate the computations.

In \cite{dennis}, the authors also considered
the branch curve of a projection of a non-prime $K3$ surface.
But they considered quotients of $\pi_1(\C\P^2-S)$
(where $S$ is the branch curve)
by a subgroup of commutators
(commutators of geometric generators
which are mapped to disjoint transpositions
by the geometric monodromy representation, see Definition 2.2 in that paper).
Their motivation came from the theory of symplectic manifolds
and families of projections
where the branch curves acquire and may lose
pairs of transverse double points with opposite orientations;
creating a pair of double points adds a commutation relation.
The quotient there is the largest quotient of $\pi_1(\C\P^2-S)$,
which is guaranteed to be invariant under these operations,
and it is different from the $\pi_1(\C\P^2-S)$
or the group related to the fundamental group of the Galois cover,
which is our aim (see \cite{pil}).

The paper is divided as follows.
In Section \ref{sec:2} we give the definition of a degeneration
and explain briefly the computations from \cite{11}.
In Section \ref{bm} we review the general setup
and the notion of the braid monodromy.
In Section \ref{branch} we recover the relevant properties
of the branch curve $S$ of the generic projection to the plane
of the $K3$ surface by using regeneration techniques.
Section \ref{Delta} states the braid monodromy factorization of $\Delta_{48}^2$ of $S$
and its related invariance properties.
We explain why the computations of $\Delta_{48}^2$
and its invariance are necessary for this work and for future work.

\section{$K3$ surfaces and their degenerations}\label{sec:2}

\subsection{$K3$ surfaces}\label{sec:2a}
Compact complex surfaces are classified into four broad categories,
based on the growth rate of sections of powers of the canonical class.
Such sections can either be always zero (rational or ruled surfaces);
form vector spaces of bounded dimension,
or have spaces of sections whose dimensions grow linearly (elliptic surfaces)
or quadratically (surfaces of general type) with the power.
For surfaces for which these sections are bounded,
some multiple of the canonical bundle is trivial,
and there are nine separate families up to complex deformation.
The surfaces of this type which are simply connected
in fact have trivial canonical bundle, and are called $K3$ surfaces;
the invariants for such surfaces are $p_g=1$, $q=0$,
$e = 24$, and $h^{1,1} = 22$, see \cite{10} and \cite{11}.
The most common example of a $K3$ surface is a smooth quartic surface in $\mathbb{P}^3$.
The moduli space of all $K3$ surfaces is $20$-dimensional.

Most $K3$ surfaces are not algebraic;
the algebraic ones are classified by an infinite collection
(depending on an integer $g \geq 2$) of $19$-dimensional moduli spaces.
The general member of the family has a rank one Picard group,
generated by an ample class $H$ with $H^2 = 2g-2$;
the general member of the linear system $|H|$ is a smooth curve of genus $g$,
and this linear system maps the $K3$ surface to $\mathbb{P}^g$
as a surface of degree $2g-2$.
For example, the quartic surfaces in $\mathbb{P}^3$
form the family with $g=3$.
The integer $g$ is called the \emph{genus} of the family.

The $K3$ surfaces may also be embedded by a multiple $cH$ of the primitive class $H$;
this will exhibit the family of $K3$ surfaces of genus $g$
as surfaces of degree $D=c^2(2g-2)$ in $\mathbb{P}^{G}$,
where $G = 1+c^2(g-1)$, whose hyperplane section
is an element of $|cH|$, and therefore is a curve of genus $G$.

\subsection{Degenerations of $K3$ surfaces}\label{sec:2b}
Let us start this section with recalling the definition of degeneration from \cite{18}.

\begin{definition} \label{df1} {\bf Projective degeneration}
Let $\Delta$ be the unit disc,
and $X, Y$ be algebraic surfaces (or more generally algebraic varieties).
Suppose that $k: Y \rightarrow \C \P^n$ and $k': X \rightarrow \C \P ^n$
are projective embeddings.
We say that $k'$ is a {\bf projective degeneration} of $k$
if there exist a flat family $\pi: V \rightarrow \Delta$,
and an embedding $F:V\rightarrow \Delta \times \C\P^n$,
such that $F$ composed with the first projection is $\pi$,
and:
\begin{itemize}
\item[(a)] $\pi^{-1}(0) \simeq X$;
\item[(b)] there is a $t_0 \neq 0$ in $\Delta$ such that
$\pi^{-1}(t_0) \simeq Y$;
\item[(c)] the family $V-\pi^{-1}(0) \rightarrow \Delta-{0}$
is smooth;
\item[(d)] restricted to $\pi^{-1}(0)$, $F = {0}\times k'$
under the identification of $\pi^{-1}(0)$ with $X$;
\item[(e)] restricted to $\pi^{-1}(t_0)$, $F = {t_0}\times k$
under the identification of $\pi^{-1}(t_0)$ with $Y$.
\end{itemize}
\end{definition}

In \cite{11}, the authors constructed a degeneration of embedded $K3$ surfaces,
to a union of $D$ planes, meeting according to the combinatorics of a certain
triangulation of the $2$-sphere.
Such degenerations are called \emph{Type III} degenerations,
and there are many such available.
The particular one constructed there was formed
by first decomposing the sphere into two rectangles,
then decomposing each rectangle into $ab$ squares
(by decomposing the sides of the rectangles into $a$ and $b$ line segments),
and finally decomposing each square into two triangles (planes).
The parameters $a$ and $b$ are free to choose,
and the construction exhibits a degeneration of $K3$ surfaces into $4ab$ planes.
If $c$ is the greatest common divisor of $a$ and $b$,
the general member of the family
is a smooth $K3$ surface embedded by the multiple $c$ of the generator $H$ of the Picard group.
Hence using the above notation, $4ab = 2c^2(g-1)$,
so that $g = 1+2(ab/c^2)$ and $G = 1+2ab$
(the number of coordinate points in $\P^G$ is $G+1$).
The boundary of the two rectangular arrays of planes contains $2a+2b$ lines,
where the identification is taking place.

This degeneration of a $K3$ surface,
is called a {\bf pillow degeneration},
see \cite{11}.

\begin{figure}[ht]
\epsfxsize=12cm 
\epsfysize=5cm 
\begin{minipage}{\textwidth}
\begin{center}
\epsfbox{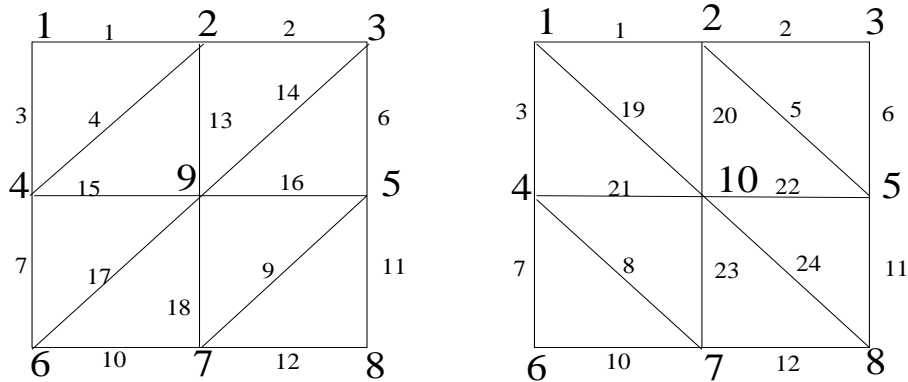}
\end{center}
\end{minipage}
\caption{The $(2,2)$-pillow degeneration}\label{top}
\end{figure}

This article proceeds by making further analyses with the $a=b=2$ case.
Given a non-prime $K3$ surface of degree $16$,
we get the pillow degeneration
$(K3)_0$, as depicted in Figure \ref{top}.
We quote the main result from \cite{11}.

\begin{theorem}
The $(2,2)$-pillow degeneration $(K3)_0$ is a total degeneration
of a smooth $K3$ surface of degree $16$ in $\P^9$.
The smooth $K3$ surface is a re-embedding of a quartic surface
in $\P^3$ via the linear system of quadrics.
\end{theorem}

For the regeneration process we have to fix a numbering of the vertices and lines.
The boundary points are labeled from $1$ to $8$,
while the interior point on the top is $9$ and on the bottom is $10$.
The $16$ planes meet each other along a total of $24$ lines,
each joining two of the $10$ coordinate points.
The numbering of the vertices
induces a lexicographic numbering of the lines as follows.
If $L$ has endpoints $a < b$ and $M$ has endpoints $c < d$
then $L < M$ if $b < d$ or $d=b$ and $a < c$.
This gives a total ordering of the lines,
which we interpret as a numbering from $1$ to $24$,
as shown in Figure \ref{top}.

Note that the coordinate points $1$, $3$, $6$, and $8$ at the corners
are each contained in three distinct planes,
while all other points are each contained in six planes.
We call these two types $3$-points and $6$-points respectively.

A general projection $f_0:(K3)_0 \to \C\P^2$ will also be a
degeneration of a general such projection of the smooth $K3$ surface.
Under $f_0$, each of the $16$ planes is mapped isomorphically to $\C\P^2$.
The ramification locus $R_0$ of $f_0$ is the closed subset of $(K3)_0$,
where $f_0$ is not a local isomorphism.
Here $R_0$ is exactly the $24$ lines.
Let $S_0 = f_0(R_0)$ be the degenerated branch curve;
it is a line arrangement,
composed of the images of the $24$ lines
(each counted twice and each of which is a line in the plane).

\bigskip
\bigskip

\section{The braid monodromy notion}\label{bm}

\begin{figure}[h]
\epsfxsize=6cm 
\epsfysize=4cm 
\begin{minipage}{\textwidth}
\begin{center}
\epsfbox{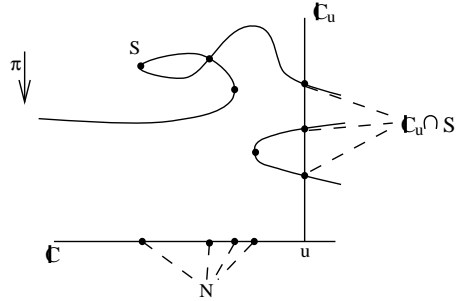}
\end{center}
\end{minipage}
\caption{General setting}\label{setup}
\end{figure}

Consider the following setting (Figure \ref{setup}).
$S$ is an algebraic curve in $\C^2$, with $p = \deg(S)$.
Let $\pi: \C^2 \rightarrow \C$ be a generic projection onto the first coordinate.
Define the fiber $K(x) = \{y \mid (x,y) \in S\}$ in $S$ over a fixed point $x$,
projected to the $y$-axis.
Define $N = \{x \mid \# K(x) < p \}$ and
$M' = \{ s \in S \mid \pi_{\mid s} \mbox{ is not \'{e}tale at } s \}$;
note that $\pi (M') = N$.
Let $\set{A_j}^q_{j=1}$ be the set of points of $M'$
and $N = \set{x_j}^q_{j=1}$ their projection on the $x$-axis.
Recall that $\pi$ is generic,
so we assume that $\# (\pi^{-1}(x) \cap M') =1$
for every $x \in N$.
Let $E$ (resp. $D$) be a closed disk on the $x$-axis (resp. the $y$-axis),
\st $M' \subset E \times D$ and $N \subset \mbox{Int}(E)$.
We choose $u \in \partial E$ a real point far enough from the set $N$,
so that $x << u$ for every $x \in N$.
Define $\C_u = \pi^{-1}(u)$
and number the points of $K=\C_u\cap S$ as $\{1 , \dots , p\}$.

We now construct a g-base for the fundamental group $\pi_1(E - N, u)$.
Take a set of paths $\set{\g_j}^q_{j=1}$
which connect $u$ with the points $\set{x_j}^q_{j=1}$ of $N$.
Now encircle each $x_j$ with a small oriented counterclockwise circle $c_j$.
Denote the path segment from $u$ to the boundary of this circle as $\g'_j$.
We define an element (a loop) in the g-base as $\delta_j = {\g'}_j c_j {\g'}^{-1}_j$.
Let $B_p[D,K]$ be the braid group,
and let $H_1, \dots , H_{p-1}$ be its frame
(for complete definitions, see \cite[Section III.2]{16}).
The braid monodromy of $S$ \cite{2} is a map
$\varphi: \pi_1(E - N, u) \rightarrow B_p[D,K]$ defined as follows:
every loop in $E - N$ starting at $u$ has liftings
to a system of $p$ paths in $(E - N) \times D$
starting at each point of $K = 1,\ldots,p$.
Projecting them to $D$ we obtain $p$ paths in $D$
defining a motion $\{1(t), \dots , p(t)\}$ (for $0 \leq t \leq 1$)
of $p$ points in $D$ starting and ending at $K$.
This motion defines a braid in $B_p [D , K]$.
By the Artin Theorem \cite{17}, for $j=1, \dots, q$,
there exists a halftwist $Z_j \in B_p[D , K]$
and $\epsilon_j \in \Z$,
\st $\varphi(\dl _j) = Z_j^{\epsilon_j}$,
where  $Z_j$ is a halftwist and $\epsilon_j = 1,2$ or $3$
(for an ordinary branch point, a node, or a cusp respectively).
We now explain how to describe this $Z_j$.

First we recall a definition of an almost real curve from \cite{17},
as we have such curves here.

\begin{definition}
A curve $S$ is called an \emph{almost real curve} if
\begin{enumerate}
\item
$N \subseteq E \cap \R$,
\item $N \subset E - \partial E$,
\item $\forall x \in {E \cap \R -N}, \# {(K \cap \R)(x)} \ge p-1$,
\item $\forall x \in N, \# {\pi^{-1}(x) \cap M'} = 1$,
\item The singularities can be
\begin{enumerate}
\item a branch point, topologically locally equivalent to $y^2+x=0$ or $y^2-x=0$,
\item a tacnode (a line is tangent to a conic)
\item an intersection of $m$ smooth branches, transversal to each other.
\end{enumerate}
\end{enumerate}
\end{definition}

We explain how to get the braid monodromy around each singularity in $S$.
Let $A_j$ be a singularity in $S$ and denote by $x_j$ its projection by $\pi$ to the $x$-axis.
We choose a point $x'_j$ next to $x_j$,
\st \ $\pi^{-1}(x'_j)$ is a typical fiber.
If $A_j$ is (b), (c) or (d),
then $x'_j$ is on the right side of $x_j$.
If $A_j$ is (a), then $x'_j$ is on the  left side of $x_j$
(the typical fiber in case (a),
which is on the left  side of this singularity,
intersects the conic in two real points).
We encircle $A_j$ with a very small circle
in such a way that the typical fiber $\pi^{-1}(x'_j)$
intersects the circle in two points, say $a ,b$.
We fix a skeleton $\xi_{x'_j}$ which connects $a$ and $b$, and denote it as $<a,b>$.
The Lefschetz diffeomorphism $\Psi$ (\cite[Subsection 1.9.5]{3})
allows us to get a resulting skeleton $(\xi_{x'_j}) \Psi$
in the typical fiber $\C_u$.
This one defines a motion of its two endpoints.
This motion induces a halftwist
$Z_j = \Delta < (\xi_{x'_j}) \Psi>$.
As above,
$\varphi(\delta_j) = \Delta < (\xi_{x'_j}) \Psi>^{\epsilon_j}$.
The braid monodromy factorization associated to $S$ is
$$
\Delta^2_{p} = \prodl^q_{j=1} \varphi(\delta_j).
$$

\section{The branch curve $S_0$ and its regeneration}\label{branch}
\def\Zovera
{\mathop{\lower 10pt \hbox{$ {\buildrel{\displaystyle \bar{Z}_{ij}}\over
{\hspace{-.2cm}\scriptstyle{(a)}}} $}}  \nolimits}

\subsection{The branch curve $S_0$}\label{ctil}
The degenerated object $(K3)_0$ has $\{j\}^{10}_{j=1}$ as vertices.
By the projection $f_0 : K3_0 \rightarrow \C\P^2$,
we obtain a line arrangement
$S_0 = \bigcup\limits^{24}_{i=1} L_i$,
and the projections $f_0(j)=j$ are singularities of $S_0$.
The $3$-points are $1, 3, 6, 8$ and the $6$-points are $2, 4, 5, 7, 9, 10$.
In the following subsection
we regenerate in neighborhoods of these singularities.
The monodromies will be given in the next section.

Besides these singularities in $S_0$,
there are also {\em parasitic intersections}:
these come from lines in $(K3)_0$ which do not intersect,
but when projecting $(K3)_0$ onto $\C\P^2$,
they will intersect.
Denote each line in $(K3)_0$ as a pair (by its two end vertices).
Let $u$ be a point, $u \notin S_0$, \st $\# \C_u = 24$,
and $q_i = i \cap \C_u$ a real point.
Take two non-intersecting lines $p = (i,k)$ and
$t = (j,\ell)$ in $(K3)_0$.
\def\Zoverthree
{\mathop{\lower 10pt \hbox{$ {\buildrel{\displaystyle\bar{Z}^2_{i \; 6}}\over
{\hspace{-.2cm}\scriptstyle{(5)}}}$}}\nolimits}
\def\Zoverfour
{\mathop{\lower 10pt \hbox{$ {
\buildrel{\displaystyle\bar{Z}^2_{i \; 5}}
\over
{\hspace{-.2cm}\scriptstyle{(4)}}} $}}  \nolimits}
\def\Zovereight
{\mathop{\lower 10pt \hbox{$ {
\buildrel{\displaystyle\bar{Z}^2_{i \; 9}}
\over{\hspace{-.2cm}\scriptstyle{(8)}}} $}}  \nolimits}
\def\Zovernine
{\mathop{\lower 10pt \hbox{$ {\buildrel{\displaystyle \bar{Z}^2_{i \; 10}}
\over{\hspace{-.2cm}\scriptstyle{(8)(9)}}} $}}  \nolimits}
\def\Zoversix
{\mathop{\lower 10pt \hbox{$ {\buildrel{\displaystyle\bar{Z}^2_{i \; 7}}
\over{\hspace{-.2cm}\scriptstyle{(6)}}} $}}  \nolimits}
\def\Zoversixseven
{\mathop{\lower 10pt \hbox{$ {\buildrel{\displaystyle\bar{Z}^2_{i \; 8}}
\over{\hspace{-.2cm}\scriptstyle{(6)(7)}}} $}}  \nolimits}
\def\Zovernineteneleven
{\mathop{\lower 10pt \hbox{$ {\buildrel{\displaystyle\bar{Z}^2_{i \; 12}}
\over{\hspace{-.2cm}\scriptstyle{(11)}}} $}}  \nolimits}
\def\Zovernineten
{\mathop{\lower 10pt \hbox{$ {\buildrel{\displaystyle\bar{Z}^2_{i \; 11}}
\over{\hspace{-.2cm}\scriptstyle{(9)(10)}}} $}}  \nolimits}
\def\Zovernineteen
{\mathop{\lower 10pt \hbox{$ {\buildrel{\displaystyle\bar{Z}^2_{i \; 20}}
\over{\hspace{-.2cm}\scriptstyle{(19)}}} $}}  \nolimits}
\def\Zoverthirteen
{\mathop{\lower 10pt \hbox{$ {\buildrel{\displaystyle\bar{Z}^2_{i \; 14}}
\over{\hspace{-.2cm}\scriptstyle{(13)}}} $}}  \nolimits}
\def\Zoverthirteenf
{\mathop{\lower 10pt \hbox{$ {\buildrel{\displaystyle\bar{Z}^2_{i \; 15}}
\over{\hspace{-.2cm}\scriptstyle{(13)(14)}}} $}}  \nolimits}
\def\Zoverthirteenff
{\mathop{\lower 10pt \hbox{$ {\buildrel{\displaystyle\bar{Z}^2_{i \; 16}}
\over{\hspace{-.2cm}\scriptstyle{(13)-(15)}}} $}}  \nolimits}
\def\Zoverseventeen
{\mathop{\lower 10pt \hbox{$ {\buildrel{\displaystyle\bar{Z}^2_{i \; 18}}
\over{\hspace{-.2cm}\scriptstyle{(13)-(17)}}} $}}  \nolimits}
\def\Zoverseventeene
{\mathop{\lower 10pt \hbox{$ {\buildrel{\displaystyle\bar{Z}^2_{i \; 19}}
\over{\hspace{-.2cm}\scriptstyle{(17)(18)}}} $}}  \nolimits}
\def\Zoverseventeenen
{\mathop{\lower 10pt \hbox{$ {\buildrel{\displaystyle\bar{Z}^2_{i \; 17}}
\over{\hspace{-.2cm}\scriptstyle{(13)-(16)}}} $}}  \nolimits}
\def\Zoverseventeentwenty
{\mathop{\lower 10pt \hbox{$ {\buildrel{\displaystyle\bar{Z}^2_{i \; 21}}
\over{\hspace{-.2cm}\scriptstyle{(19)(20)}}} $}}  \nolimits}
\def\Zovercheers
{\mathop{\lower 10pt \hbox{$ {\buildrel{\displaystyle\bar{Z}^2_{i \; 24}}
\over{\hspace{-.2cm}\scriptstyle{(19)-(23)}}} $}}  \nolimits}
\def\Zoverseventeentwt
{\mathop{\lower 10pt \hbox{$ {\buildrel{\displaystyle\bar{Z}^2_{i \; 22}}
\over{\hspace{-.2cm}\scriptstyle{(19)-(21)}}} $}}  \nolimits}
\def\Zovertwentytwothree
{\mathop{\lower 10pt \hbox{$ {\buildrel{\displaystyle\bar{Z}^2_{i \; 23}}
\over{\hspace{-.2cm}\scriptstyle{(19)-(22)}}} $}}  \nolimits}
\def\Zovertwothreefour
{\mathop{\lower 10pt \hbox{$ {\buildrel{\displaystyle\bar{Z}^2_{i \; 25}}
\over{\hspace{-.2cm}\scriptstyle{(22)-(24)}}} $}}  \nolimits}
\def\Zovertwofive
{\mathop{\lower 10pt \hbox{$ {\buildrel{\displaystyle\bar{Z}^2_{i \; 26}}
\over{\hspace{-.2cm}\scriptstyle{(22)-(25)}}} $}}  \nolimits}
\def\Zovertwosix
{\mathop{\lower 10pt \hbox{$ {\buildrel{\displaystyle\bar{Z}^2_{i \; 27}}
\over{\hspace{-.2cm}\scriptstyle{(22)-(26)}}} $}}  \nolimits}
\def\Zoverxyz
{\mathop{\lower 10pt \hbox{$ {\buildrel{\displaystyle\bar{Z}^2_{i \; 9}}
\over{\hspace{-.2cm}\scriptstyle{(8)}}} $}}  \nolimits}

\begin{notation}
We denote by $Z_{ij}$ (resp. $\bZ_{ij}$)
the counterclockwise halftwist of $i$ and $j$
along a path below (resp. above) the real axis.
If $P$ is a set of points between $i$ and $j$,
$\stackrel{(P)}{Z_{ij}}$ denotes the path from $i$ to $j$
going above the points in $P$ and below the points not in $P$.
Conjugation of braids is defined as $a^b = b^{-1}ab$.
\end{notation}

We first compute the products
$D_t : = \prodl_{\stackrel{p<t}{p \cap t = \emptyset}} \tilde{Z}^2_{pt}$
as explained in \cite[Theorem IX.2.1]{16}.

\begin{align*}
\ D_1 =& D_2 = D_4 = Id \ , \ D_3 = {Z}^2_{2 \; 3} \ , \ \ D_5 =  \bZ^2_{3 \; 5} \ , \
\ D_6 =  \prodl_{i=1,3,4} \Zoverthree \ , \
\ D_7 =  \prodl_{i=1,2,5,6} \;\; \  \bar{Z}^2_{i \; 7}\\
\ D_8 =& \prodl_{i=1,2,5,6} \bZ^2_{i \; 8} \ , \
\ D_9  = \prodl_{i=1-4,7} \; \Zoverxyz \ , \
\ D_{10}=   \prodl_{i=1}^6 \; \Zovernine \ , \
\ D_{11}=   \prodl_{i=1-4,7,8,10} \;\; \bar{Z}^2_{i \; 11} \ , \\
\ D_{12} =&  \prodl^7_{i=1} \;\; \Zovernineteneleven \ , \
\ D_{13} = \prodl^{12}_{\stackrel{\scriptstyle i=3}{\scriptstyle i \neq 4,5}} \;
            \bar{Z}^2_{i \; 13} \ , \
\ D_{14} =  \prodl^{12}_{\stackrel{\scriptstyle i=1}{\scriptstyle i \neq 2,6}} \;\;
            \Zoverthirteen \ , \
\ D_{15} =  \prodl^{12}_{\stackrel{\scriptstyle i=1}{\scriptstyle i \neq 3,4,7,8}} \;\;
            \Zoverthirteenf \ , \\
\ D_{16} =&  \prodl^{12}_{\stackrel{\scriptstyle i=1}{\scriptstyle i \neq 5,6,9,11}} \;\;
            \Zoverthirteenff \ , \
\ D_{17} = \prodl^{12}_{\stackrel{\scriptstyle i=1}{\scriptstyle i \neq 7,10}} \;\;
            \Zoverseventeenen \ , \
\ D_{18} =   \prodl^{11}_{\stackrel{\scriptstyle i=1}{\scriptstyle i \neq 8-10}} \;\;
             \Zoverseventeen \ , \
\ D_{19} =  \prodl^{18}_{\stackrel{\scriptstyle i=2}{\scriptstyle i \neq 3}} \;
            \bar{Z}^2_{i \; 19} \ , \\
\ D_{20} =&  \prodl^{18}_{\stackrel{\scriptstyle i=3}{\scriptstyle i \neq 4,5,13}} \;\;
            \Zovernineteen \ , \
\ D_{21} = \prodl^{18}_{\stackrel{\scriptstyle i=1}{\scriptstyle i \neq 3,4,7,8,15}} \;\;
            \Zoverseventeentwenty \ , \
\ D_{22} =  \prodl^{18}_{\stackrel{\scriptstyle i=1}{\scriptstyle i \neq 5,6,9,11,16}} \;\;
            \Zoverseventeentwt  \ , \\
\ D_{23} =&  \prodl^{17}_{\stackrel{\scriptstyle i=1}{\scriptstyle i \neq 8-10,12}} \;\;
            \Zovertwentytwothree \ , \
\ D_{24} =  \prodl^{18}_{\stackrel{\scriptstyle i=1}{\scriptstyle i \neq 11,12}} \;\;
            \Zovercheers.
\end{align*}

Define the parasitic intersection braids as
\begin{equation}\label{formulac}
\tilde{C}_j = \prodl_{j \in t}D_t,
\end{equation}
where  $j$ is the \underline{smallest} endpoint in the line $t$.
In our case,
$\tC_1 = D_1 \cdot D_3  \cdot D_{19}, \ \ \tC_2 = D_2 \cdot D_4  \cdot D_{5}
\cdot D_{13} \cdot D_{20}, \ \ \tC_3 = D_6 \cdot D_{14}, \ \ \tC_4 = D_7 \cdot D_{8}
\cdot  D_{15} \cdot D_{21}, \ \ \tC_5 =  D_{9} \cdot  D_{11} \cdot  D_{16} \cdot  D_{22} ,
\ \ \tC_6 =  D_{10} \cdot  D_{17}, \ \ \tC_7 =  D_{12} \cdot  D_{18} \cdot D_{23},
\ \ \tC_8 =  D_{24}, \ \ \tC_9 = \tC_{10} =Id.$

\subsection{The regeneration of $S_0$}\label{regcurve}
The degenerated branch curve $S_0$ of $(K3)_0$ has degree $24$.
However each of the $24$ lines of $S_0$ should be counted as
a double line in the scheme-theoretic branch locus,
since it arises from a line of nodes.
Another way to see this is to note that
the regeneration of $(K3)_0$ induces a regeneration of $S_0$
in such a way that each point, say $c$, on the typical fiber
is replaced by two nearby points $c, c'$.

The curve $S_0$ has $3$-points, $6$-points and parasitic intersections.
In the forthcoming subsections we explain how to regenerate the curve
in neighbourhoods of these points.
The resulting branch curve $S$ will have degree $48$.

\subsubsection{Regeneration of $3$-points}\label{3pt}
The $3$-points in $S_0$ are $1, 3, 6, 8$, see Figure \ref{top}.
The regeneration is divided into steps.
We explain each step in two levels,
first dealing with the surface and then with the branch curve.

In the surface level,
each diagonal is replaced with a conic by a partial regeneration.
Focusing on a $3$-point,
we have a partial regeneration of two of the planes to a quadric surface.
We get one quadric and one plane, which is tangent to the quadric.
The plane and the quadric meet along two lines
(one from each ruling of the quadric).

For the regeneration of the branch curve,
we need the following lemma from \cite{19}.

\begin{lemma}\label{lm:24}
Let $V$ be a projective algebraic surface, and $D'$ be a curve in $V$.
Let $f: V \rightarrow \C \P^2$ be a generic projection.
Let $S \subseteq \C \P^2, S' \subset V$
be the branch curve of $f$ and the corresponding ramification curve.
Assume $S'$ intersects $D'$ at a point $\alpha'$.
Let $D = f(D')$ and $\alpha = f(\alpha ')$.
Assume that there exist neighbourhoods of $\alpha $ and $\alpha '$,
\st $f_{\mid_{S'}}$ and $f_{\mid_{D'}}$ are isomorphisms.
Then $D$ is tangent to $S$ at $\alpha $.
\end{lemma}

At the branch curve level,
we have two double lines
(coming from the intersection of the plane and the quadric)
and one conic (coming from the branching of the quadric over the plane).
According to the above lemma,
the conic is tangent to each of the two double lines.

\begin{figure}[ht]
\epsfxsize=10cm 
\epsfysize=4cm 
\begin{minipage}{\textwidth}
\begin{center}
\epsfbox{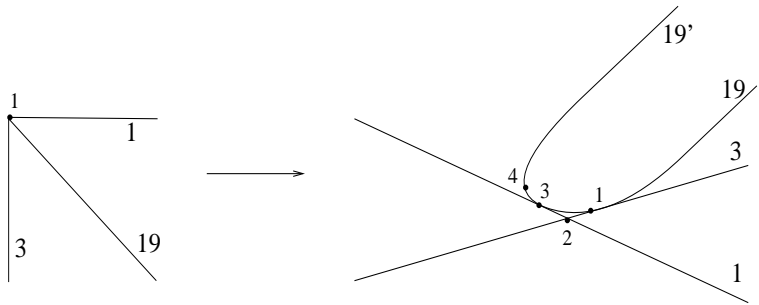}
\end{center}
\end{minipage}
\caption{Partial regeneration in a neighbourhood of the $3$-point $1$}\label{1reg1}
\end{figure}

As far as the branch points go,
one of the two branch points of the conic is far away from the $3$-point,
and the other one is close to the $3$-point; see for example
Figure \ref{1reg1}.

In the next step of the regeneration, at the curve level,
we use regeneration Lemmas from \cite{17}.
The two tangent points regenerate to three cusps each (giving a total of six)
and the intersection point of the two double lines gives eight more branch points.
One can think of this as first giving four nodes,
then each node giving two branch points.

At the surface level it means that we get a smooth surface
which locally looks like a cubic in $\C\P^3$
(degenerating to a triple of planes).

\subsubsection{Regeneration of $6$-points}\label{6pt}
The $6$-points in $S_0$ are $2, 4, 5, 7, 9, 10$.
Each of these points is the projection to the plane
of a cone over a cycle of independent lines spanning a $\C\P^5$.
Hence, locally, the surface lives in $\C\P^6$ and consists of six planes through a point.

\begin{figure}[ht]
\epsfxsize=10cm 
\epsfysize=5cm 
\begin{minipage}{\textwidth}
\begin{center}
\epsfbox{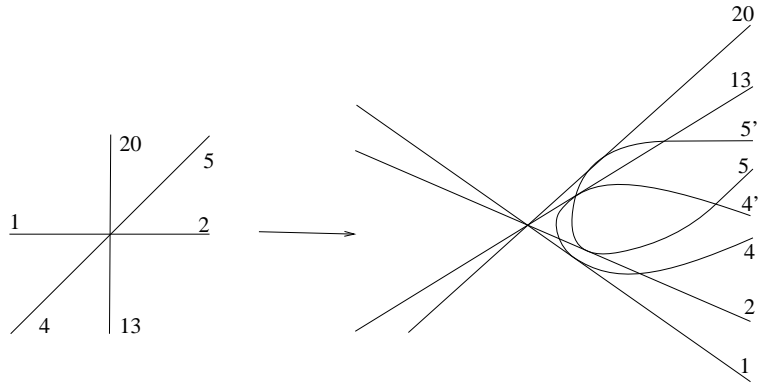}
\end{center}
\end{minipage}
\caption{Partial regeneration in a neighbourhood of the $6$-point $2$}\label{2phi}
\end{figure}

Let us understand how the regeneration of a $6$-point occurs,
by focusing for example on a neighbourhood around $2$ in $S_0$.
The first regeneration consists in smoothing out the lines $4$ and $5$ and
therefore replacing four of these pairwise adjacent but opposite planes with two
quadrics, see Figure \ref{2phi}.

\begin{figure}
\epsfxsize=7cm 
\epsfysize=3.5cm 
\begin{minipage}{\textwidth}
\begin{center}
\epsfbox{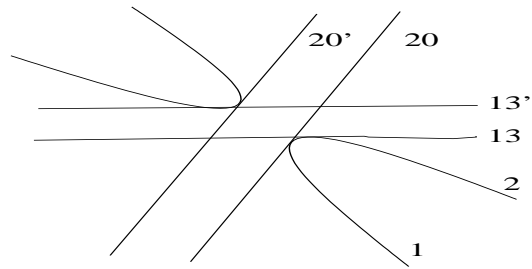}
\end{center}
\end{minipage}
\caption{Second step of regeneration}\label{4lines}
\end{figure}

Then one smooths out the lines $13$ and $20$.
This produces two cubic rational normal scrolls meeting along the lines $1$ and $2$.
Now the union of these lines on each scroll is homologous to a non singular conic,
so we can further regenerate to make the two scrolls meet along a smooth conic
(in this way $1$ and $2$ are replaced by a smooth double conic).
We visualize a purely local figure (Figure \ref{4lines}),
in order to understand this step.

Finally one smooths the conic too and arrives to a Del-Pezzo sextic in $\C\P^6$.

A similar regeneration occurs in the neighbourhoods of the other $6$-points.

\section{The braid monodromy factorization  $\Delta^2_{48}$}\label{Delta}
In Section \ref{branch} we explained how to obtain the branch curve $S$,
which has degree $48$.
At this point we will compute the braid monodromy factorization $\Delta^2_{48}$ of $S$.

In Subsection \ref{rules} we formulate braids and regeneration rules.
In Subsection \ref{delta} we compute $\Delta^2_{48}$.
In Subsections \ref{prop} and \ref{imp}
we show an invariance property of $\Delta^2_{48}$,
and emphasise its importance.

\subsection{The regeneration rules}\label{rules}
Recall that $Z_{ij}$ is a counterclockwise halftwist of two points $i$ and $j$.
We start with an example that illustrates conjugated braids.

\begin{example}\label{example6}
The path on the left-hand side in Figure \ref{example} is constructed as follows:
take a path $z_{34}$ and conjugate it by the fulltwist $\bZ_{13}^2$
(1 encircles 3 counterclockwise while moving above the axis).
We get the left-hand side path $z_{34}^{\bZ_{13}^2}$.
The right-hand side path is constructed as follows:
take again $z_{34}$ and conjugate it first by $Z_{23}^2$
(3 encircles 2 counterclockwise) and then by $Z_{13}^2$
(3 encircles 1 counterclockwise while moving below the axis).
We get the right-hand side path $z_{34}^{Z_{23}^2Z_{13}^2}$.
The related halftwists of $z_{34}^{\bZ_{13}^2}$ and $z_{34}^{Z_{23}^2Z_{13}^2}$
are $Z_{34}^{\bZ_{13}^2}$ and $Z_{34}^{Z_{23}^2Z_{13}^2}$ respectively.

\begin{figure}[h!]
\epsfxsize=7cm 
\epsfysize=1cm 
\begin{minipage}{\textwidth}
\begin{center}
\epsfbox{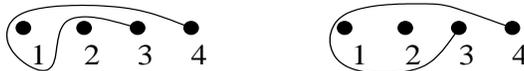}
\end{center}
\end{minipage}
\caption{Example of conjugated braids}\label{example}
\end{figure}
\end{example}

In Section \ref{regcurve} we explained that
by the regeneration one gets the original branch curve $S$.
Let $N, M, u$ be as in Section \ref{bm}.
The set $K = S \cap \C_u$ are the intersection points of the curve $S$
with the typical fiber $\C_u$; $\# K = 48$.
Recall that by the regeneration,
each point $c$ in $K_0 = S_0 \cap \C_u$ is replaced by two close points in $K$,
say $c, c'$.

We define the braid $Z^2_{ij}$ to be a fulltwist of $j$ around $i$,
and $Z^2_{i'j}$ to be a fulltwist of $j$ around $i'$.
The braid $Z^2_{ii',j}$ is obtained by a regeneration
(the point $i$ on the typical fiber is replaced by $i, i'$)
and it is a fulltwist of $j$ around $i$ and $i'$.
This braid (and similar ones) is formulated in the following lemma:

\begin{lemma}{\bf \cite[Lemma 2.5]{3}}\\
The following formulas hold:  $Z^2_{ii',j} =  Z^2_{i'j} Z^2_{ij},
\ \ Z^{2}_{i',jj'}  =  Z^{2}_{i'j'} Z^{2}_{i'j}, \ \ Z^{-2}_{i',jj'}  =
Z^{-2}_{i'j} Z^{-2} _{i'j'}, \ \ {\bar{Z}}^{-2}_{i',jj'}  =  {\bar{Z}}^{-2}_{i'j'}
{\bar{Z}}^{-2}_{i'j}, \ \ Z^{-2}_{ii',j}  =  Z^{-2}_{ij} Z^{-2} _{i'j}, \ \
Z^2_{ii',jj'} = Z^2_{i',jj'} Z^2_{i,jj'}, \ \
Z^{-2}_{ii',jj'}  = Z^{-2}_{i,jj'} Z^{-2}_{i',jj'}$.
\end{lemma}

In the following regeneration rules we shall describe what happens to a
factor in a factorized expression of a braid monodromy by a regeneration.

\begin{theorem}\label{1rule}
{\bf First regeneration rule \cite[p. 336]{19}}\\
A factor of the form $Z_{ij}$ regenerates to $Z_{ij'} Z_{i'j}$.
\end{theorem}

\begin{theorem}\label{2rule}
{\bf Second regeneration rule \cite[p. 337]{19}}\\
A factor of the form $Z^2_{ij}$ regenerates to $Z^2_{ii',j}, Z^2_{i,jj'}$
or $Z^2_{ii',jj'}$.
\end{theorem}

\begin{theorem}\label{3rule}
{\bf Third regeneration rule \cite[p. 337]{19}}\\
A factor of the form $Z^4_{ij}$ regenerates to
$Z^{3}_{i,jj'} = (Z^3_{ij})^{Z_{jj'}} \cdot (Z^3_{ij}) \cdot (Z^3_{ij})^{{Z_{jj'}}^{-1}}$
or to $Z^{3}_{ii',j}=(Z^3_{ij})^{Z_{ii'}} \cdot (Z^3_{ij}) \cdot (Z^3_{ij})^{{Z_{ii'}}^{-1}}$.
\end{theorem}

\subsection{The factorization}\label{delta}
The braid monodromy factorization $\Delta^2_{48}$ of $S$ is
$\prodl^{10}_{i=1} C_i \varphi_{i}$,
where $C_i$ are the regenerations
of the parasitic intersection braids from Section \ref{ctil},
and $\varphi_{i}$ are the local braid monodromies
which we get when regenerating around the singularities $1, \dots, 10$.

\subsubsection{Braid monodromies related to $3$-points}\label{mono3}
The $3$-points are $1, 3, 6, 8$.
We concentrate in the neighbourhood of $1$,
see Figure \ref{1reg1}.
First the diagonal line $19$ regenerates to a conic $Q_{19}$
which is tangent to the two other lines $1$ and $3$, see Lemma \ref{lm:24}.
We compute the braid monodromy $\varphi_1$ of the resulting curve.

\begin{proposition}\label{14}
The local braid monodromy $\varphi_1$ is
\begin{align*}
\varphi_1 = & Z^3_{3 \; 3', 19} \cdot {(Z_{1\ 3} \cdot Z_{1\ 3} \cdot  Z_{1'\ 3} \cdot
Z_{1'\ 3} \cdot Z_{1\ 3'} \cdot Z_{1\ 3'} \cdot  Z_{1'\ 3'} \cdot Z_{1'\ 3'})}^{Z^2_{3 \; 3', 19}}
\cdot \\
& {Z}^3_{1 \; 1', 19}  \cdot {Z_{19\ 19'}}^{Z^2_{3 \; 3', 19} {Z}^2_{1 \; 1', 19}}.
\end{align*}
\end{proposition}

\begin{proof}
We follow Figure \ref{1reg1}.
Let $\pi_1 : E \times D \rightarrow E$ be the projection to $E$.

Let $\{j\}^{4}_{j=1}$ be singular points of $\pi_1$ as follows:
$1$, $3$ are the tangent points of $Q_{19}$ with the lines $L_3$, $L_1$ respectively,
$2$ is the intersection  point of the  lines $L_1$, $L_3$, and
$4$ is the branch point in $Q_{19}$.

Let $N = \{x(j) = x_j \mid 1 \le j \le 4 \}$,
\st $N \subset E - \partial E, N \subset E$.
Take $u \in \partial E$, \st $\C_u$ is a typical fiber and $M \in \C_u$ is a real point.
Recall that $K = K(M)$.
$K = \{1, 3, 19, 19'\}$, \st the points are  real and
$1 < 3 < 19 < 19'$.
Let $i = L_i \cap K$ for $i = 1, 3$ and
$\{19, 19'\} = Q_{19} \cap K$.

We are looking for $\vp_M(\delta_j)$ for $j = 1,\ldots,4$.
So we choose a $g$-base $\{\delta_j\}^{4}_{j=1}$ of $\pi_1(E - N,u)$, \st
each $\delta_j$ is constructed from a path $\g_j$ below the real line
and a counterclockwise small circle around the points in $N$.

The diffeomorphism which is induced from passing through a branch point
was defined in \cite{3} and \cite{17} by  $\Delta^{\frac{1}{2}}_{I_2\R}<k>$.
We recall the precise definition from \cite{3}.
Consider a typical fiber on the left side of a branch point
(locally defined by $y^2-x=0$).
The typical fiber intersects the conic in two complex points.
Passing through this point,
these points become real on the right-hand side typical fiber.
The two points  move to the k'th place and rotate in a counterclockwise $90^o$ twist.
They become real and numbered as $k, k+1$.

First find the skeleton $\xi_{x'_j}$ related to each singular point, as
explained in Section \ref{bm}.
Then compute the local diffeomorphisms $\delta_{j}$
induced from singular points $j$.
Since the points $1$ and $3$ (resp. $2$) are
tangent points (resp. a node), the diffeomorphisms $\delta_{1}$ and $\delta_{3}$
are each of degree $2$ (resp. $1$).
\vspace{-1cm}
\begin{center}
\begin{tabular}{cccc} \\
$j$ & $\xi_{x'_j}$ & $\epsilon_{j}
$ & $\delta_{j}$ \\ \hline
1 & $<3,19>$ & 4 & $\Delta^2<3,19>$\\
2 & $<1,3>$ & 2 & $\Delta<1,3>$\\
3 & $<3,19>$ & 4 & $\Delta^2<3,19>$\\
4 & $<19,19'>$ & 1 & $\Delta^{\frac{1}{2}}_{I_2\R}<19>$\\
\end{tabular}
\end{center}

Using \cite[Theorems 1.41, 1.44]{3} and \cite{17},
we compute the skeleton  $(\xi_{x'_j}) \Psi_{\g'_j}$ to each $j$
by applying to the skeleton $\xi_{x'_j}$
the product $\prodl_{i=j-1}^{1} \delta_{i}$.
\begin{list}{}{\leftmargin 2cm \labelsep .4cm \labelwidth 3cm}
\item $(\xi_{x'_1}) \Psi_{\g'_1} =<3,19>= z_{3\ 19}$\\
$\varphi_M(\delta_1) = Z^4_{3\ 19}$
\vspace{-.3cm}
\begin{figure}[h!]\label{1phiAA}
\epsfysize=0.5cm 
\begin{minipage}{\textwidth}
\begin{center}
\epsfbox{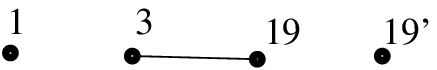}
\end{center}
\end{minipage}
\end{figure}
\item $(\xi_{x'_2}) \Psi_{\g'_2} = <1,3>
\Delta^2<3,19>= z_{1\ 3 }^{Z^2_{3\ 19}}$\\
$\varphi_M(\delta_2) = {Z^2_{1\ 3 }}^{Z^2_{3\ 19}}$
\begin{figure}[h!]\label{1phiBB}
\epsfysize=1.5cm 
\begin{minipage}{\textwidth}
\begin{center}
\epsfbox{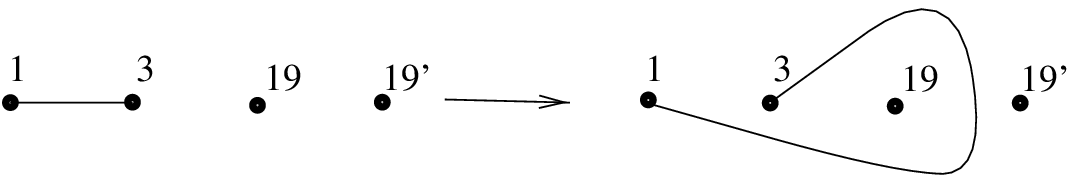}
\end{center}
\end{minipage}
\end{figure}
\item $(\xi_{x'_3}) \Psi_{\g'_3}=<3,19>\Delta<1,3>\Delta^2<3,19>={z}_{1\ 19}$\\
$\varphi_M(\delta_3) = {Z}^4_{1\ 19}$
\begin{figure}[h!]\label{1phiCC}
\epsfysize=1.3cm 
\begin{minipage}{\textwidth}
\begin{center}
\epsfbox{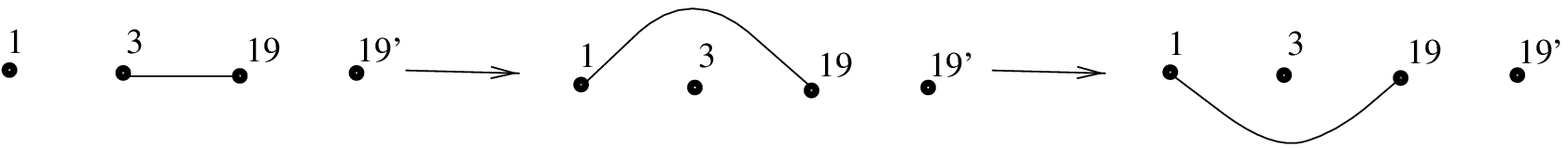}
\end{center}
\end{minipage}
\end{figure}
\item $(\xi_{x'_4}) \Psi_{\g'_4} =
<19,19'>\Delta^2<3,19>\Delta<1,3>\Delta^2<3,19>=z_{19\ 19'}^
{Z^2_{3\ 19} {Z}^2_{1\ 19}}$\\
$\varphi_M(\delta_4) = {Z_{19\ 19'}}^{Z^2_{3\ 19} {Z}^2_{1\ 19}}$
\begin{figure}[ht]\label{1phiDD}
\epsfysize=1.5cm 
\begin{minipage}{\textwidth}
\begin{center}
\epsfbox{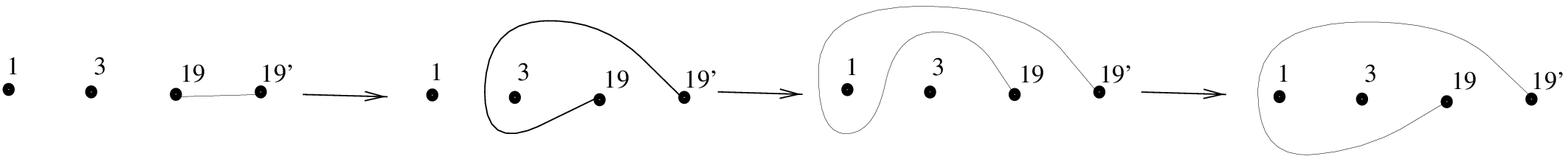}
\end{center}
\end{minipage}
\end{figure}
\end{list}

Now, by Theorem \ref{3rule},
in the regeneration process each one of the tangent points
regenerates to three cusps.
Therefore, the factors $Z^4_{3\ 19}$ and ${Z}^4_{1\ 19}$
regenerate to $Z^3_{3 \; 3', 19}$ and ${Z}^3_{1 \; 1', 19}$ respectively.
As explained in Subsection \ref{3pt},
the node is replaced by eight branch points:
since the lines $1$ and $3$ double to be $1$, $1'$ and $3$, $3'$ respectively,
we get four nodes, the intersections of $1$ and $3$,  $1'$ and $3$,  $1$ and $3'$,
$1'$ and $3'$.
Each node regenerates to two branch points.
For $1$ and $3$ we get easily the braid monodromy $Z_{1\ 3} \cdot Z_{1\ 3}$.
For the other ones we get $Z_{1'\ 3} \cdot Z_{1'\ 3}$,
$Z_{1\ 3'} \cdot Z_{1\ 3'}$ and $Z_{1'\ 3'} \cdot Z_{1'\ 3'}$.
In Figure \ref{1phi} we show the paths which are related to the braids in $\varphi_1$.
\begin{figure}[h!]
\epsfxsize=14.5cm 
\epsfysize=6cm 
\begin{minipage}{\textwidth}
\begin{center}
\epsfbox{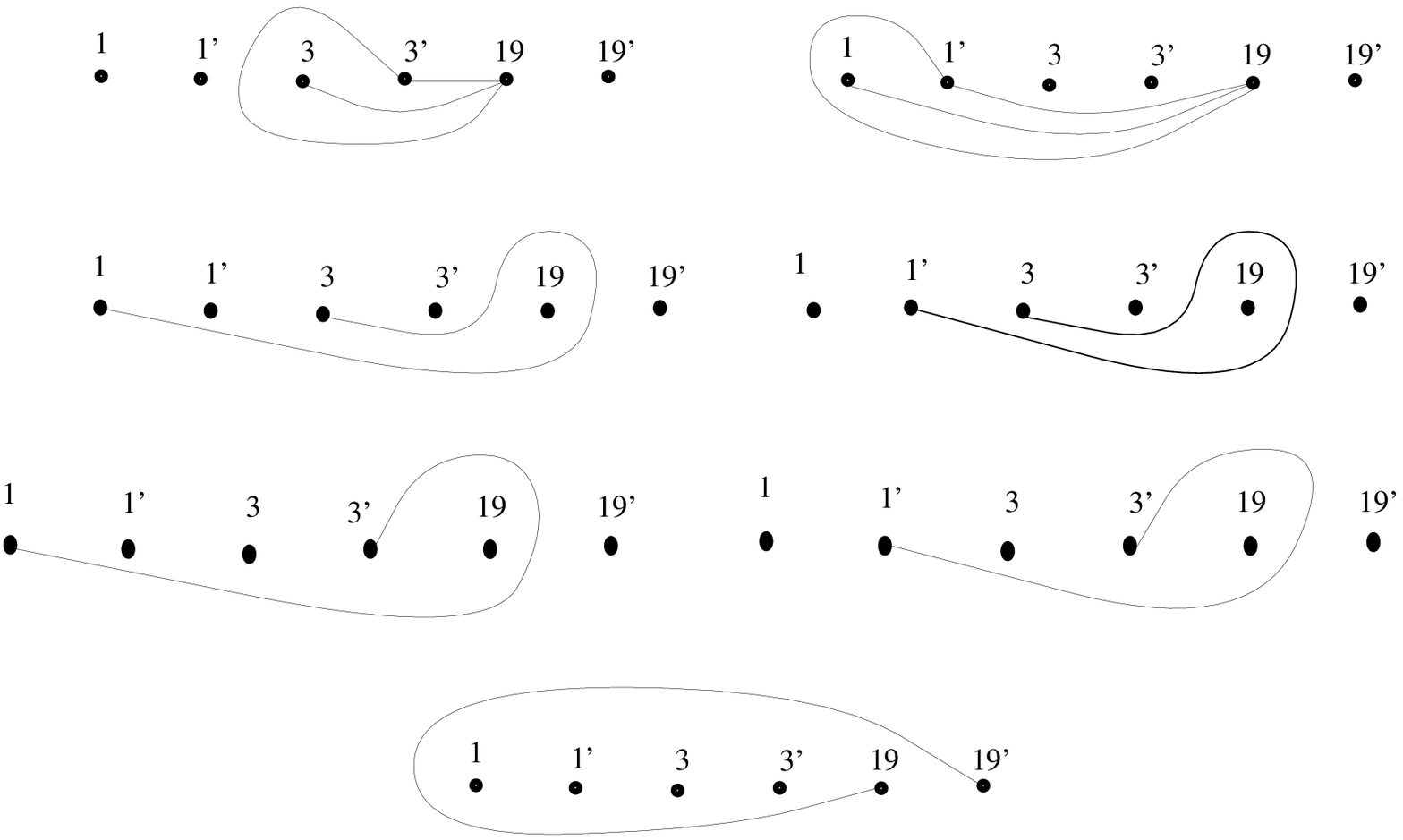}
\end{center}
\end{minipage}
\caption{}\label{1phi}
\end{figure}

\end{proof}

The computations for the points $3, 6$ and $8$ are the same.
One has to exchange the indices $1, 3$ and $19$ to $2, 6$ and $14$ (for the point $3$),
to $7, 10$ and $17$ (for the point $6$)
and to $11, 12$ and $24$ (for the point $8$) respectively.
Therefore we can formulate their related monodromies.
\begin{theorem}\label{second2}
The local braid monodromy around the $3$-point $m = 1, 3, 6$ or $8$ is formulated as follows:
\begin{eqnarray*}
{}\varphi_m =& Z^3_{j \; j', k} \cdot {(Z_{i\ j} \cdot Z_{i\ j} \cdot Z_{i'\ j} \cdot
Z_{i'\ j} \cdot Z_{i\ j'} \cdot Z_{i\ j'} \cdot Z_{i'\ j'} \cdot Z_{i'\ j'})}^{Z^2_{j \; j', k}}
\cdot {Z}^3_{i \; i', k} \cdot {Z_{k\ k'}}^{Z^2_{j \; j', k} Z^2_{i \; i', k}},
\end{eqnarray*}
where $i, j$ and $k$ are the three lines which meet at the $3$-point $m$ and satisfy $i < j < k$.
\end{theorem}

\subsubsection{Braid monodromies related to $6$-points}\label{mono6}
The $6$-points in $S_0$ are $2, 4, 5, 7, 9, 10$.
In Section \ref{6pt} we recalled the regeneration near the $6$-point $2$.
The regenerations near $4, 5, 7, 9, 10$ are similar to the one around $2$,
but they differ in the indices.

Regenerations in the neighbourhood of a $6$-points
were studied carefully in \cite{3}, \cite{5}, \cite{8} and \cite{19}.
Therefore, we only state the resulting monodromies.

\begin{theorem}\label{mono2}
The local braid monodromy $\varphi_2$ has the following form
\begin{align*}
\varphi_2 = & Z^2_{2 \; 2',4} \cdot Z^2_{4',5 \; 5'}  \cdot Z^3_{1 \;1',4} \cdot
{Z^3_{4',13 \;13'}}^{Z^{2}_{4' ,5 \; 5'}} \cdot {Z^2_{2 \;2',4}}^{Z^{2}_{1 \; 1',4}} \cdot
{Z^2_{4',5 \;5'}}^{Z^{-2}_{5 \; 5',13\; 13'}} \cdot \\
& {Z_{4 \; 4'}}^{Z^2_{4',13 \; 13'}Z^2_{4',5 \;5'}Z^2_{2 \; 2',4}Z^2_{1 \; 1',4}} \cdot
{Z^2_{4',20\; 20'}}^{Z^{2}_{4',13\; 13'}Z^{2}_{4',5\; 5'}} \cdot \\
& {Z^2_{4,20\; 20'}}^{Z^{2}_{4,13\; 13'}Z^{2}_{4,5\; 5'}Z^2_{4 \; 4'}Z^2_{2 \;2',4}
Z^{2}_{1 \; 1',4}}
\cdot \left(G \cdot {\left(F \cdot F^{Z_{20 \; 20'}^{-1}Z_{13 \; 13'}^{-1}}\right)}^
{Z^2_{2 \; 2',5}Z^2_{1 \; 1',5}}\right)^{Z^2_{2 \; 2',4}Z^2_{1 \; 1',4}},
\end{align*}
where
\begin{align*}
G = & Z^2_{5',13 \; 13'}  \cdot Z^3_{2 \;2',5} \cdot {Z^3_{5',20 \;20'}}^{Z^{2}_{5' ,13 \; 13'}}
\cdot {Z^2_{5',13 \;13'}}^{Z^{-2}_{13 \; 13' ,20 \; 20'}} \cdot
{Z_{5 \; 5'}}^{Z^2_{2 \; 2',5}Z^2_{5',20 \;20'}Z^2_{5',13 \; 13'}} \cdot \\
& Z^2_{1 \; 1',5} \cdot {Z^2_{1 \; 1',5'}}^{Z^2_{5',20 \; 20'}Z^2_{5',13 \;13'}}
\end{align*}
and
$$
F = Z^3_{2 \; 2',13} \cdot Z^2_{13' \; 20}\cdot {Z^2_{13 \; 20}}^{Z^2_{2 \; 2',13}}\cdot
{Z^3_{2 \; 2',20}}^{Z^2_{2 \; 2',13}} \cdot \tilde{Z}_{1 \; 2'} \cdot \tilde{Z}_{1' \; 2}.
$$
The figures which correspond to the braids outside $G$ and $F$ are \ref{v2a}-\ref{v2i}.
The ones which correspond to the factors in $G$ and $F$ are \ref{v2j}-\ref{v2p} and
\ref{f2a}-\ref{f2e1} respectively.
\begin{figure}[ht]
\begin{minipage}{\textwidth}
\begin{center}
\epsfbox{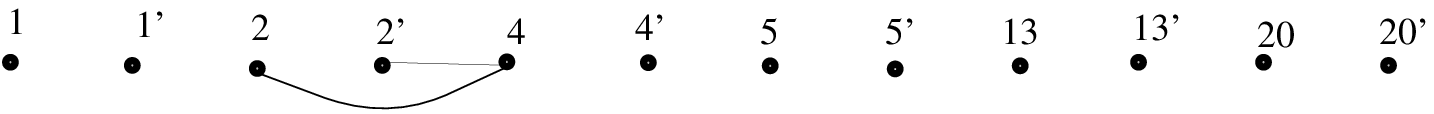}
\end{center}
\end{minipage}
\caption{}\label{v2a}
\end{figure}
\begin{figure}[ht]
\begin{minipage}{\textwidth}
\begin{center}
\epsfbox{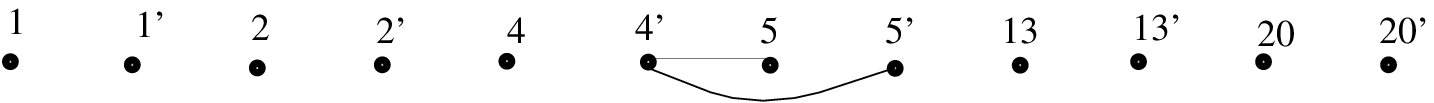}
\end{center}
\end{minipage}
\caption{}\label{v2b}
\end{figure}
\begin{figure}[ht]
\begin{minipage}{\textwidth}
\begin{center}
\epsfbox{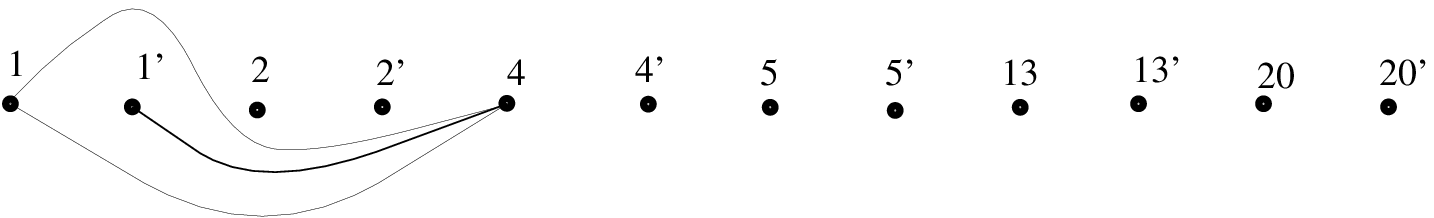}
\end{center}
\end{minipage}
\caption{}\label{v2c}
\end{figure}
\begin{figure}[ht]
\begin{minipage}{\textwidth}
\begin{center}
\epsfbox{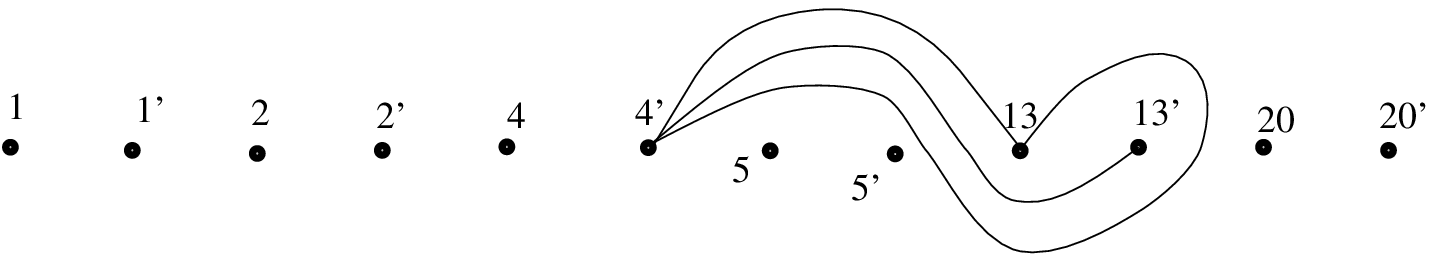}
\end{center}
\end{minipage}
\caption{}\label{v2d}
\end{figure}
\begin{figure}[ht]
\begin{minipage}{\textwidth}
\begin{center}
\epsfbox{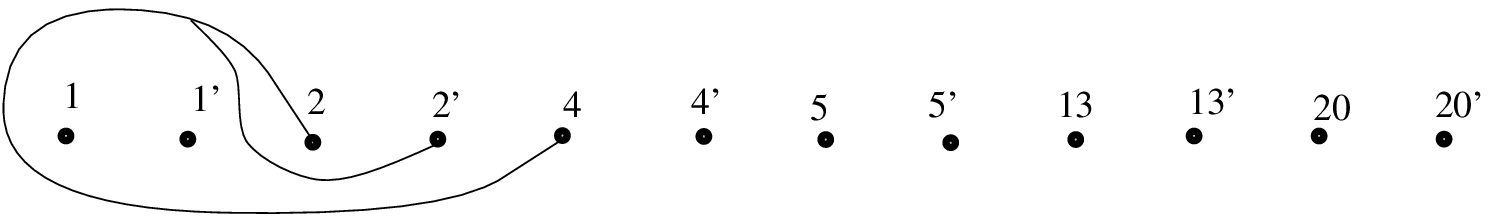}
\end{center}
\end{minipage}
\caption{}\label{v2e}
\end{figure}
\begin{figure}[ht]
\begin{minipage}{\textwidth}
\begin{center}
\epsfbox{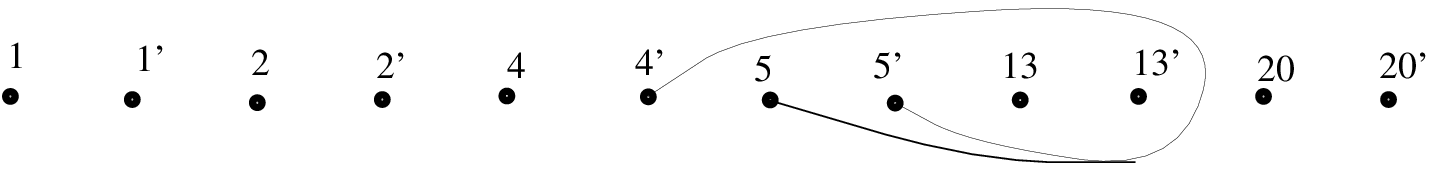}
\end{center}
\end{minipage}
\caption{}\label{v2f}
\end{figure}
\begin{figure}[ht]
\begin{minipage}{\textwidth}
\begin{center}
\epsfbox{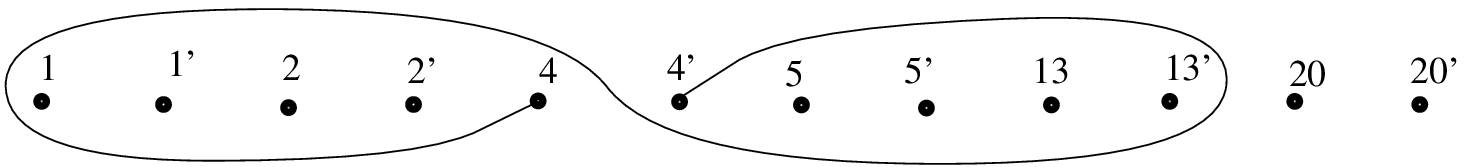}
\end{center}
\end{minipage}
\caption{}\label{v2g}
\end{figure}
\begin{figure}[ht]
\begin{minipage}{\textwidth}
\begin{center}
\epsfbox{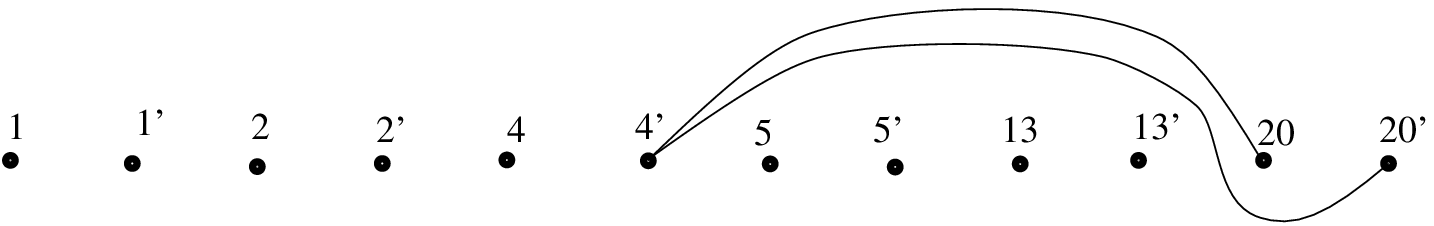}
\end{center}
\end{minipage}
\caption{}\label{v2h}
\end{figure}
\begin{figure}[ht]
\begin{minipage}{\textwidth}
\begin{center}
\epsfbox{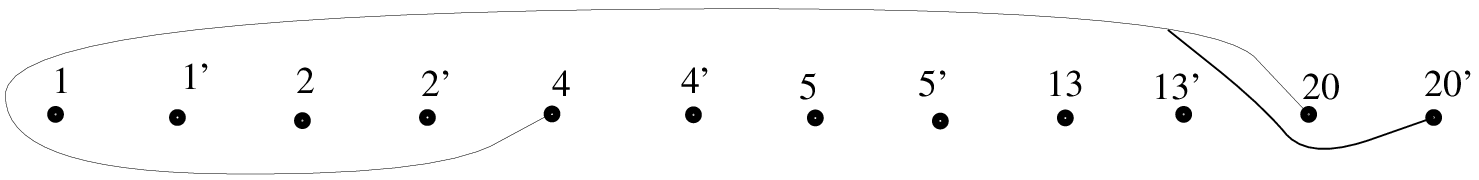}
\end{center}
\end{minipage}
\caption{}\label{v2i}
\end{figure}
\begin{figure}[ht]
\begin{minipage}{\textwidth}
\begin{center}
\epsfbox{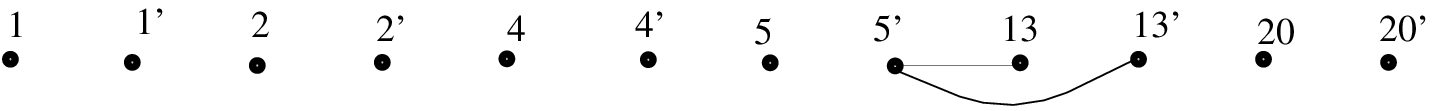}
\end{center}
\end{minipage}
\caption{}\label{v2j}
\end{figure}
\begin{figure}[ht]
\begin{minipage}{\textwidth}
\begin{center}
\epsfbox{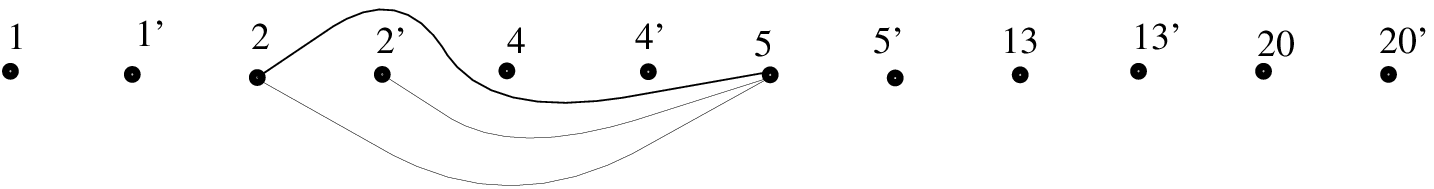}
\end{center}
\end{minipage}
\caption{}\label{v2k}
\end{figure}
\begin{figure}[ht]
\begin{minipage}{\textwidth}
\begin{center}
\epsfbox{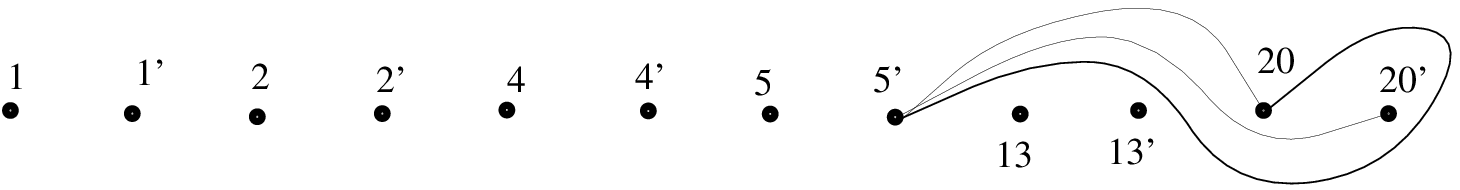}
\end{center}
\end{minipage}
\caption{}\label{v2l}
\end{figure}
\begin{figure}[ht]
\begin{minipage}{\textwidth}
\begin{center}
\epsfbox{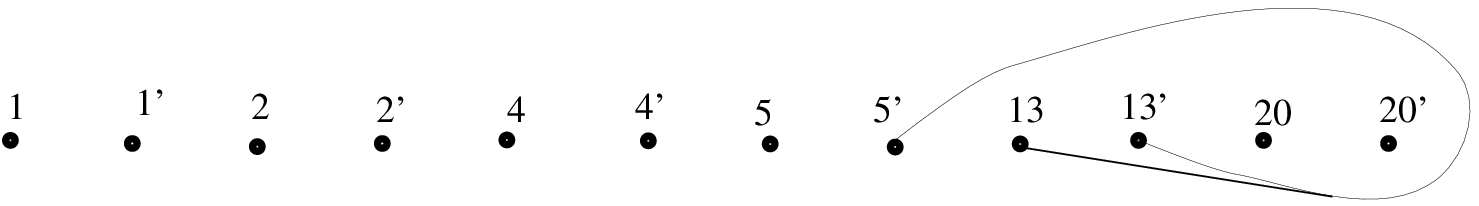}
\end{center}
\end{minipage}
\caption{}\label{v2m}
\end{figure}
\begin{figure}[ht]
\begin{minipage}{\textwidth}
\begin{center}
\epsfbox{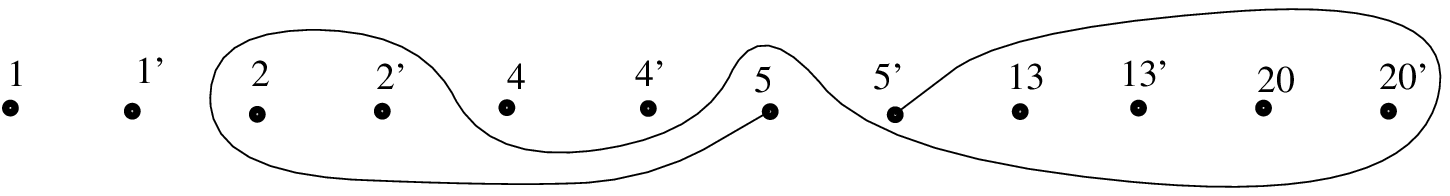}
\end{center}
\end{minipage}
\caption{}\label{v2n}
\end{figure}
\begin{figure}[ht]
\begin{minipage}{\textwidth}
\begin{center}
\epsfbox{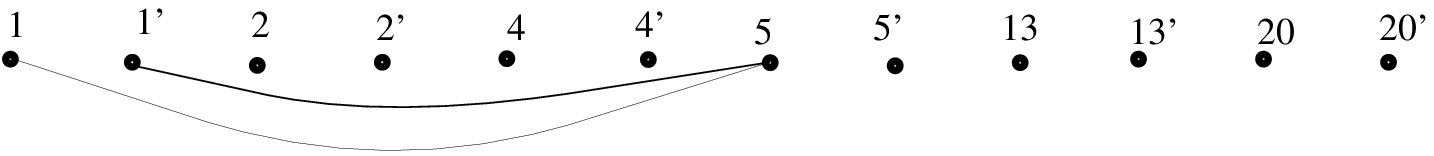}
\end{center}
\end{minipage}
\caption{}\label{v2o}
\end{figure}
\begin{figure}[ht]
\begin{minipage}{\textwidth}
\begin{center}
\epsfbox{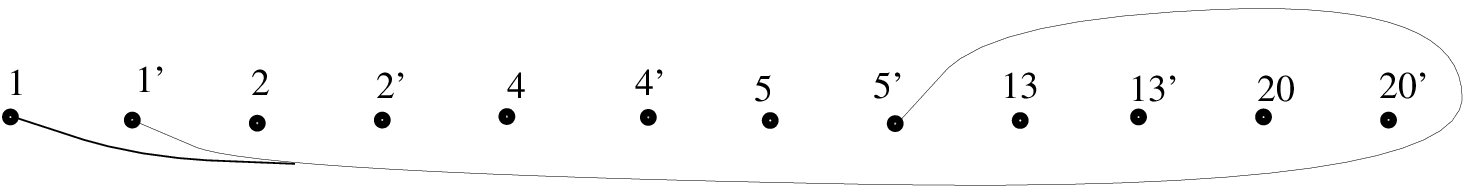}
\end{center}
\end{minipage}
\caption{}\label{v2p}
\end{figure}
\begin{figure}
\begin{minipage}{\textwidth}
\begin{center}
\epsfbox{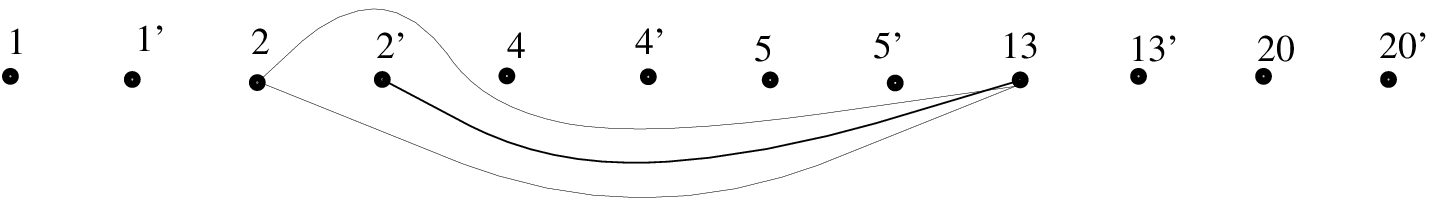}
\end{center}
\end{minipage}
\caption{}\label{f2a}
\end{figure}
\begin{figure}
\begin{minipage}{\textwidth}
\begin{center}
\epsfbox{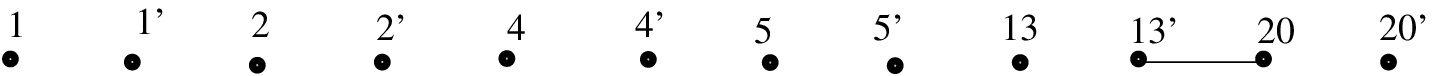}
\end{center}
\end{minipage}
\caption{}\label{f2c1}
\end{figure}
\clearpage
\begin{figure}
\begin{minipage}{\textwidth}
\begin{center}
\epsfbox{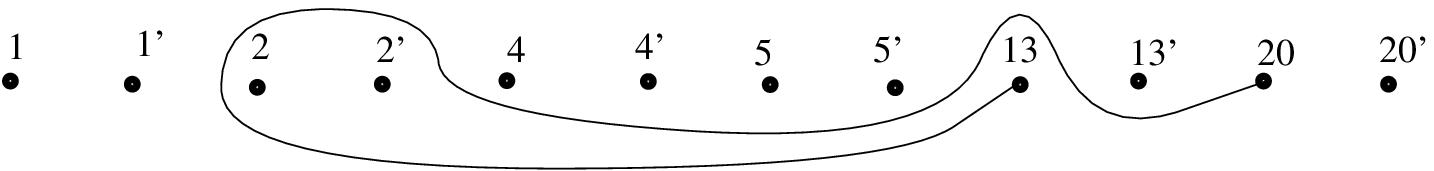}
\end{center}
\end{minipage}
\caption{}\label{f2c2}
\end{figure}
\begin{figure}
\begin{minipage}{\textwidth}
\begin{center}
\epsfbox{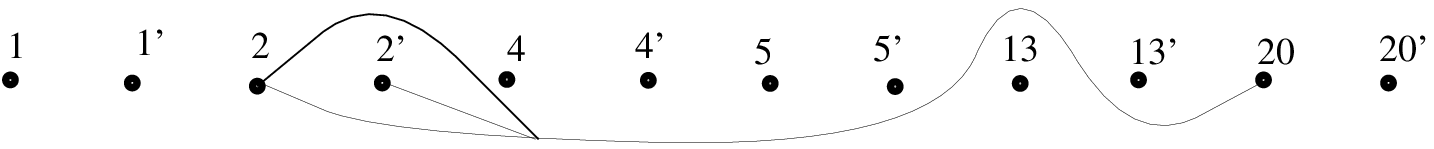}
\end{center}
\end{minipage}
\caption{}\label{f2e}
\end{figure}
\begin{figure}
\begin{minipage}{\textwidth}
\begin{center}
\epsfbox{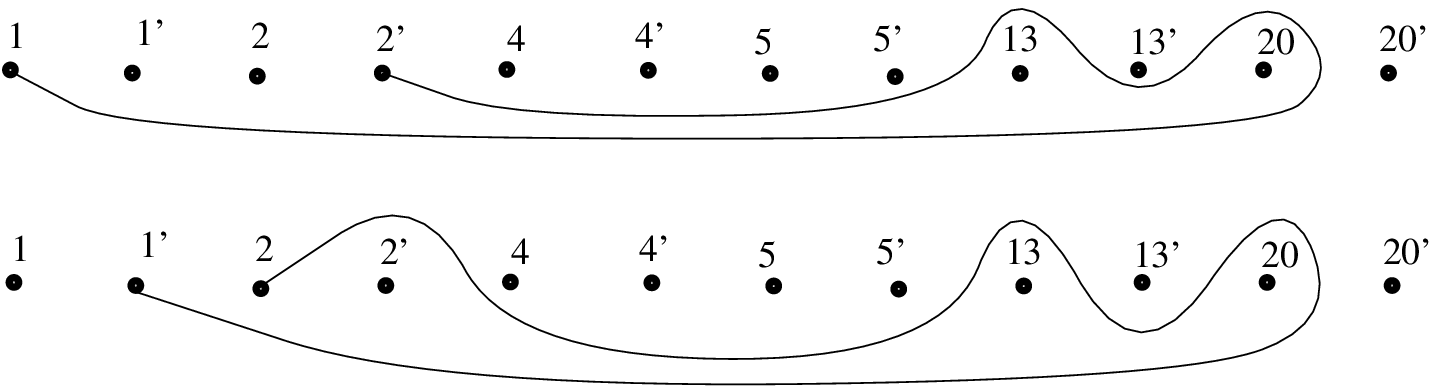}
\end{center}
\end{minipage}
\caption{}\label{f2e1}
\end{figure}

\end{theorem}

In the following theorem we give the local braid monodromies
$\varphi_4, \varphi_5, \varphi_7, \varphi_9$ and $\varphi_{10}$.
Related figures can be easily constructed, using Example \ref{example6}.
We note that $\varphi_{10}$ was computed precisely in \cite{3} and \cite{19}
(and related figures can be found there).

\begin{theorem}
\begin{enumerate}
\item
The local braid monodromy $\varphi_4$ has the following form
\begin{align*}
\varphi_4 = & Z^2_{4',7 \; 7'}  \cdot {Z^2_{4',8 \; 8'}}^{Z^{2}_{4' ,7 \; 7'}} \cdot
Z^3_{3 \; 3',4} \cdot {Z^3_{4',15 \; 15'}}^{Z^{2}_{4',8\; 8'}Z^{2}_{4',7\; 7'}} \cdot
{Z^2_{4',8 \; 8'}}^{Z^{2}_{4',7\; 7'}Z^{-2}_{8\; 8',15\; 15'}} \cdot \\
& {Z^2_{4',7 \; 7'}}^{Z^{-2}_{7\; 7', 8\; 8'}Z^{-2}_{7\; 7',15\; 15'}} \cdot
{Z_{4 \; 4'}}^{Z^2_{3 \; 3',4}Z^2_{4',15 \; 15'}Z^2_{4',8\; 8'}Z^2_{4',7 \;7'}} \cdot
{Z^2_{4,21\; 21'}}^{Z^{2}_{4,15\; 15'}Z^{2}_{4,8\; 8'}Z^2_{4,7 \; 7'}Z^2_{4 \;4'}Z^{2}_{3 \; 3',4}}
\cdot \\
& {Z^2_{4',21\; 21'}}^{Z^{2}_{4',15\; 15'}Z^{2}_{4',8\; 8'}Z^2_{4',7 \; 7'}} \cdot
\left(G \cdot {\left(F \cdot F^{Z_{7 \; 7'}^{-1}Z_{3 \; 3'}^{-1}}\right)}^
{Z^2_{7 \; 7',8}Z^2_{3 \; 3',8}}\right)^{Z^2_{3 \; 3',4}},
\end{align*}
where
\begin{align*}
G = & Z^2_{8',15 \; 15'}  \cdot Z^3_{7 \;7',8} \cdot {Z^3_{8',21 \;21'}}^{Z^{2}_{8',15 \; 15'}}
\cdot {Z^2_{8',15 \; 15'}}^{Z^{-2}_{15 \; 15' ,21 \; 21'}} \cdot
{Z_{8 \; 8'}}^{Z^2_{7 \; 7',8}Z^2_{8',21 \;21'}Z^2_{8',15 \; 15'}} \cdot \\
& Z^2_{3 \; 3',8} \cdot {Z^2_{3 \; 3',8'}}^{Z^2_{8',21 \; 21'}Z^2_{8',15 \;15'}}
\end{align*}
and
$$
F=Z^2_{3' \; 7} \cdot Z^3_{7',15 \; 15'} \cdot {Z^2_{3' \; 7'}}^
{Z^2_{7',15 \; 15'}Z^2_{3' \; 7}} \cdot {Z^3_{3',15 \; 15'}}^{Z^2_{3' \; 7}} \;
\cdot {\tilde{Z}_{ 15 \; 21'}} \cdot {\tilde{Z}_{ 15' \; 21}}.
$$
The paths related to ${\tilde{Z}_{ 15 \; 21'}}$ and ${\tilde{Z}_{ 15' \; 21}}$
appear in Figure \ref{alpha4}.
\begin{figure}[ht]
\begin{minipage}{\textwidth}
\begin{center}
\epsfbox{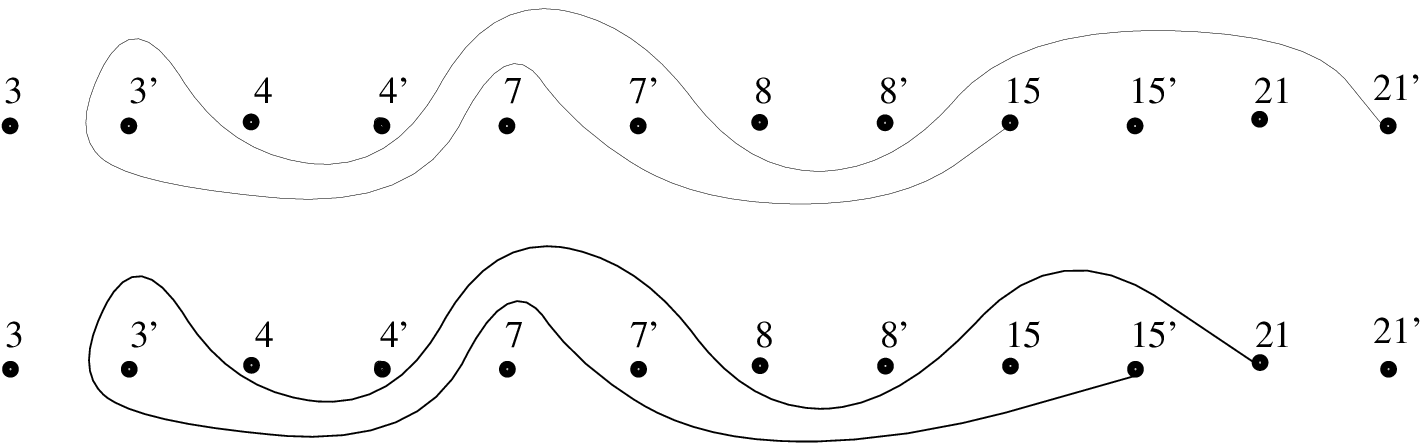}
\end{center}
\end{minipage}
\caption{}\label{alpha4}
\end{figure}
\item
The local braid monodromy $\varphi_5$ has the following form
\begin{align*}
\varphi_5 = & Z^3_{5',6 \; 6'}  \cdot G^{Z^2_{5',6 \; 6'}} \cdot
{Z^3_{5',22 \; 22'}}^{Z^{2}_{5',6\; 6'}} \cdot
{Z_{5 \; 5'}}^{Z^{2}_{5',22\; 22'}Z^{2}_{5',6\; 6'}} \cdot
{Z^2_{5',16 \; 16'}}^{Z^{2}_{5\; 5'}Z^{2}_{5',6\; 6'}} \cdot \\
& {Z^2_{5',16 \; 16'}}^{Z^2_{5',6 \; 6'}} \cdot
{Z^2_{5,11\; 11'}}^{Z^{2}_{5\; 5'}Z^{2}_{5',6\; 6'}} \cdot
{Z^2_{5',11\; 11'}}^{Z^{2}_{5',6\; 6'}} \cdot
{Z^2_{5,9\; 9'}}^{Z^{2}_{5\; 5'}Z^{2}_{5',6\; 6'}} \cdot
{Z^2_{5',9\; 9'}}^{Z^{2}_{5',6\; 6'}},
\end{align*}
where
\begin{align*}
G = & Z^2_{6\; 6',9}  \cdot Z^3_{9', 11\;11'} \cdot {Z^2_{6\; 6',9'}}^{Z^{2}_{6 \; 6',9}Z^{2}_{9',11
\; 11'}}  \cdot {\left(F \cdot F^{Z_{11 \; 11'}^{-1}
Z_{6 \; 6'}^{-1}}\right)}^{Z^2_{6 \; 6',9 \; 9'}Z^2_{9', 11\; 11'}} \cdot \\
& {Z^2_{9',22 \; 22'}}^{Z^{2}_{9',11 \; 11'}} \cdot
Z^2_{9, 22 \; 22'} \cdot {Z^3_{9', 16 \; 16'}}^{Z^2_{9',11 \; 11'}} \cdot
{Z_{9 \; 9'}}^{Z^{2}_{9',16\; 16'}Z^{2}_{9',11\; 11'}}
\end{align*}
and
$$
F=  Z^2_{6' \; 11} \cdot Z^3_{11',16 \; 16'} \cdot {Z^2_{6' \; 11'}}^
{Z^2_{11',16 \; 16'}Z^2_{6' \; 11}} \cdot {Z^3_{6',16 \; 16'}}^{Z^2_{6' \; 11}}
\cdot \tilde{Z}_{ 16 \; 22'} \cdot \tilde{Z}_{ 16' \; 22}.
$$
The paths related to $\tilde{Z}_{ 16 \; 22'}$ and $\tilde{Z}_{ 16' \; 22}$
are the same ones as in Figure \ref{alpha4}, with the following exchange of indices:
$3 \mapsto 6, 7 \mapsto 11, 15 \mapsto 16, 21 \mapsto 22$.
\item
The local braid monodromy $\varphi_7$ has the following form
\begin{align*}
\varphi_7 = & Z^2_{8 \; 8',9}  \cdot Z^2_{9', 10\; 10'} \cdot
{Z^3_{9',12 \; 12'}}^{Z^{2}_{9',10\; 10'}} \cdot
{Z^2_{9',10 \; 10'}}^{Z^{-2}_{10\; 10',12 \; 12'}} \cdot
{Z^2_{8\; 8',9'}}^{Z^{2}_{9',12\; 12'}Z^{2}_{9',10\; 10'}Z^{2}_{8\; 8',9}} \cdot  \\
& G^{Z^2_{8\; 8',9 \; 9'}Z^2_{9', 12\; 12'}Z^2_{9', 10\; 10'}} \cdot
{Z^2_{9',23 \; 23'}}^{Z^{2}_{9',12\; 12'}Z^{2}_{9',10\; 10'}} \cdot
{Z^3_{9',18\; 18'}}^{Z^{2}_{9',12\; 12'}Z^{2}_{9',10\; 10'}} \cdot Z^2_{9,23\; 23'} \cdot \\
& {Z_{9\; 9'}}^{Z^{2}_{9',18\; 18'}Z^2_{9',12\; 12'}Z^{2}_{9',10\; 10'}},
\end{align*}
where
\begin{align*}
G = & Z^3_{8', 10\; 10'}  \cdot  {\left(F \cdot F^{Z_{23 \; 23'}^{-1}Z_{18 \; 18'}^{-1}}\right)}^{Z^2_{8',10 \; 10'}}
 \cdot {Z^3_{8', 23\;23'}}^{Z^{2}_{8',10 \; 10'}}  \cdot
{Z_{8 \; 8'}}^{Z^{2}_{8',23\; 23'}Z^{2}_{8',10\; 10'}} \cdot \\
& {Z^2_{8',18 \; 18'}}^{Z^{2}_{8',10 \; 10'}} \cdot {Z^2_{8, 18 \; 18'}}^{Z^2_{8 \; 8'}Z^2_{8',10 \; 10'}} \cdot
{Z^2_{8', 12 \; 12'}}^{Z^2_{8',10 \; 10'}} \cdot
{Z^2_{8, 12 \; 12'}}^{Z^2_{8 \; 8'}Z^2_{8',10 \; 10'}}
\end{align*}
and
$$
F = Z^3_{12 \; 12',18} \cdot Z^2_{18' \; 23} \cdot {Z^2_{18 \; 23}}^{Z^2_{12 \; 12',18}} \cdot
{Z^3_{12 \; 12',23}}^{Z^2_{12 \; 12',18}} \cdot \tilde{Z}_{10 \; 12'} \cdot \tilde{Z}_{10' \; 12}.
$$
Moreover, the paths related to $\tilde{Z}_{10 \; 12'}$ and $\tilde{Z}_{10' \; 12}$
are the same ones as in Figure \ref{f2e1}, with the following exchange of indices:
$1 \mapsto 10, 2 \mapsto 12, 13 \mapsto 18, 20 \mapsto 23$.
\item
The local braid monodromy $\varphi_9$ has the following form
\begin{eqnarray*}
\varphi_9 = & Z^2_{14',15 \; 15'}  \cdot Z^3_{13 \; 13',14} \cdot
{\bar{Z}}^3_{14', 16 \; 16'} \cdot {Z^2_{14',15
\;15'}}^{Z^{-2}_{15 \; 15',16\; 16'}} \cdot {Z_{14 \;
14'}}^{Z^{2}_{13\; 13',14}Z^{2}_{14',16\; 16'}Z^{2}_{14',15\;
15'}} \cdot \\ & \left ( {\bar{Z}}^2_{14,17 \;17'}\right
)^{Z^2_{13\;13',14}} \cdot {\bar{Z}}^2_{14',17 \; 17'} \cdot
\left({\bar{Z}^2}_{14,18\;18'}\right)^{Z^2_{13\; 13',14}}
\cdot {\bar{Z}}^2_{14',18\; 18'} \cdot \\
& \left(G \cdot {\left(F \cdot F^{Z_{18 \; 18'}^{-1}Z_{13 \; 13'}^{-1}}\right)}^
{Z^{-2}_{17, 18 \; 18'}}\right)^{Z^2_{13 \; 13',14}},
\end{eqnarray*}
where
\begin{eqnarray*}
G = & Z^2_{16 \; 16',17}  \cdot Z^3_{17',18 \;18'} \cdot {Z^3_{15 \; 15',17}} \cdot
{Z^2_{16 \; 16',17}}^{Z^{2}_{15 \; 15' ,17}} \cdot
{Z_{17 \; 17'}}^{Z^2_{17',18 \; 18'}Z^2_{16 \;16',17}Z^2_{15 \; 15',17}} \cdot \\
& Z^2_{13 \; 13',17} \cdot {Z^2_{13 \; 13',17'}}^{Z^2_{17',18 \; 18'}}
\end{eqnarray*}
and
$$
F=Z^3_{13',15 \; 15'} \cdot Z^3_{16 \; 16',18} \cdot  \tilde{Z}_{15 \; 16'} \cdot
\tilde{Z}_{15' \; 16} \cdot {Z^2_{13' \; 18}}^{Z^2_{13',15 \; 15'}} \cdot Z^2_{13 \; 18}.
$$
The paths related to $\tilde{Z}_{15 \; 16'}$ and $\tilde{Z}_{15' \; 16}$
appear in Figure \ref{alpha9}.
\begin{figure}[ht!]
\begin{minipage}{\textwidth}
\begin{center}
\epsfbox{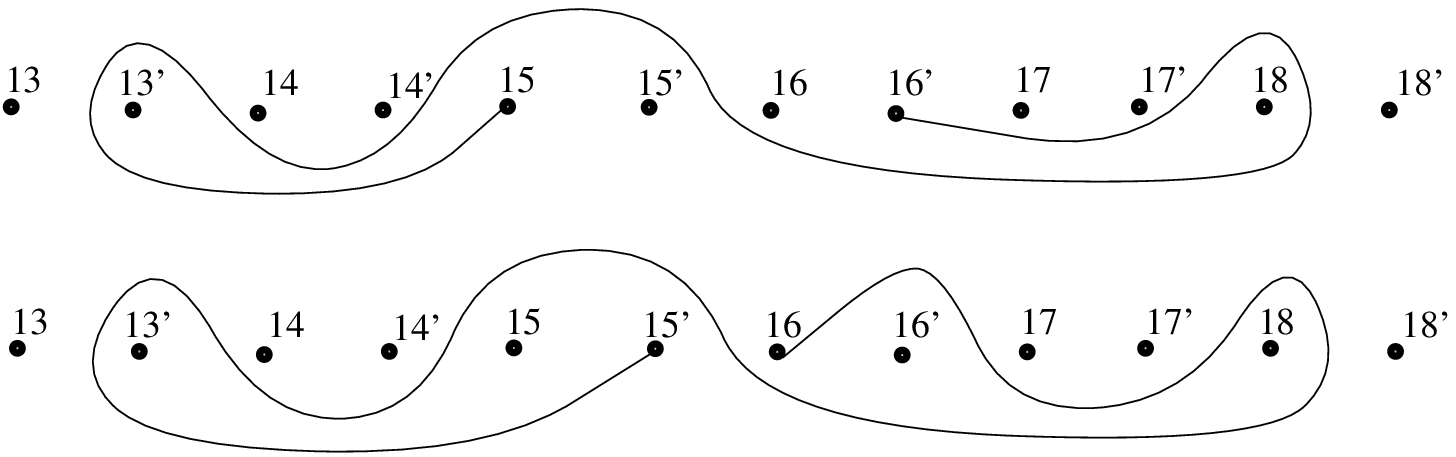}
\end{center}
\end{minipage}
\caption{}\label{alpha9}
\end{figure}
\item
The local braid monodromy $\varphi_{10}$ has the following form
\begin{align*}
\varphi_{10}= & Z^3_{19',20\; 20'} \cdot {Z_{24\; 24'}}^{Z_{22\; 22',24}^{-2}Z_{23\; 23',24}^{-2}} \cdot
Z^2_{21\; 21',24'} \cdot {\bar{Z}}^3_{22\; 22',24} \cdot
{Z^2_{21\; 21',24}}^{Z^{-2}_{23\; 23',24}} \cdot \\
& G^{Z^{2}_{19', 20 \; 20'}} \cdot Z^3_{23\; 23',24}  \cdot {\bar{Z}}^3_{19',21\; 21'} \cdot
Z^2_{19, 22\; 22'} \cdot Z^2_{19, 23\; 23'} \cdot Z^2_{19, 24\; 24'},
\end{align*}
where
\begin{align*}
G = & Z^2_{20\; 20',24'} \cdot {Z^2_{20\; 20',24}}^{Z^{-2}_{23\; 23',24}}
\cdot \left(F \cdot F^{Z_{23 \; 23'}^{-1}Z_{20 \; 20'}^{-1}}\right)^{Z^{-2}_{23 \; 23',24}}  \cdot
 Z^2_{19', 22\; 22'}  \cdot \\
& Z^2_{19', 23\; 23'} \cdot Z^2_{19', 24\; 24'}  \cdot {Z_{19\; 19'}}^{Z^2_{19',21\; 21'}}
\end{align*}
and
$$
F = Z^3_{20',21 \; 21'} \cdot Z^3_{22 \; 22',23} \cdot \tilde{Z}_{21\; 22'} \cdot \tilde{Z}_{21'\; 22}
\cdot {Z^2_{20' \; 23}}^{Z^2_{20',21 \; 21'}} \cdot Z^2_{20 \; 23}.
$$
The paths related to $\tilde{Z}_{21\; 22'}$ and $\tilde{Z}_{21'\; 22}$
are the same ones as in Figure \ref{alpha9}, with the following exchange of indices:
$13 \mapsto 20, 15 \mapsto 21, 16 \mapsto 22, 18 \mapsto 23$.
\end{enumerate}
\end{theorem}

\subsubsection{The products $C_i$}\label{regci}
The deformation also regenerates the parasitic intersection braids $\set{\tC_i}_{i = 1}^{10}$,
see Theorem \ref{2rule}.
By \cite[Lemma 19]{19}, we can take the complex conjugations
\begin{align*}
D_1 =&  D_2 = D_4 = Id \ ,
\ D_3 = Z^2_{2 \; 2',3 \; 3'} \ ,
\ D_5 = {Z}^2_{3\; 3',5 \; 5'} \ ,
\ D_{6} = \prod\limits_{i=1,3,4}^{} \stackrel{\scriptstyle (5)(5')}{Z^2}_
{\hspace{-.4cm}i\; i',6 \; 6'} \ , \\
\ D_{7} =& \prod\limits_{i=1,2,5,6}^{} {Z^2}_{\hspace{-.2cm}i\; i',7 \; 7'} \ ,
\ D_{8} = \prod\limits_{i=1,2,5,6}^{} Z^2_{i\; i',8 \; 8'} \ ,
\ D_{9} = \prod\limits_{i=1-4,7}^{} \stackrel{\scriptstyle (8)(8')}{Z^2}_
{\hspace{-.4cm}i\; i',9 \; 9'} \ ,
\ D_{10} = \prod\limits_{i=1}^6 \stackrel{\scriptstyle (8)-(9')}{Z^2}_
{\hspace{-.4cm}i\; i',10 \; 10'} \ , \\
\ D_{11} =& \prod\limits_{\stackrel{i=1}{i\neq 5,6,9}}^{10}{Z^2}_{\hspace{-.1cm}i\; i',11 \; 11'} \ ,
\ D_{12} = \prod\limits^7_{i=1} \stackrel{\scriptstyle (11)(11')}
{Z^2}_{\hspace{-.5cm}i\; i',12 \; 12'} \ ,
\ D_{13} = \prod\limits^{12}_{\stackrel{i=3}{i \neq 4,5}}
Z^2_{i\; i', 13 \; 13'} \ ,\\
\ D_{14} =& \prod\limits^{12}_{\stackrel{i=1} {i \neq 2,6}}
\stackrel{\scriptstyle (13)(13')}{Z^2}_{\hspace{-.4cm}i\; i',14 \; 14'} \ ,
\ D_{15} = \prod\limits^{12}_{\stackrel{i=1}{i \neq 3,4,7,8}}
\stackrel{\scriptstyle (13)-(14')}{Z^2}_{\hspace{-.6cm}i\; i',15 \; 15'} \ ,
\ D_{16} = \prod\limits^{12}_ {\stackrel{i=1}{i \neq 5,6,9,11}} \stackrel
{\scriptstyle (13)-(15')} {Z^2}_{\hspace{-.6cm}i\; i',16 \; 16'} \ ,\\
\ D_{17} =& \prod\limits^{12}_{\stackrel{i=1}{i \neq 7,10}}
\stackrel{\scriptstyle (13)-(16')}{Z^2}_{\hspace{-.6cm}i\; i',17 \; 17'} \ ,
\ D_{18} = \prod\limits^{11}_{\stackrel{i=1}{i \neq 8-10}}\stackrel{\scriptstyle (13)-(17')}
{Z^2}_{\hspace{-.6cm}i\; i',18 \; 18'} \ ,
\ D_{19} = \prod\limits^{18}_{\stackrel{i=2}{i \neq 3}} Z^2_{i\; i', 19 \; 19'} \ ,\\
\ D_{20} =& \prod\limits^{18}_{\stackrel{i=3}{i \neq 4,5,13}}
\stackrel{\scriptstyle (19)(19')}{Z^2}_{\hspace{-.4cm}i\; i',20 \; 20'} \ ,
\ D_{21} = \prod\limits^{18}_{\stackrel{i=1}{i \neq 3,4,7,8,15}}
\stackrel{\scriptstyle (19)-(20')} {Z^2}_{\hspace{-.6cm}i\; i',21 \; 21'} \ ,
\ D_{22} = \prod\limits^{18}_{\stackrel{i=1}{i \neq 5,6,9,11,16}}
\stackrel{\scriptstyle (19)-(21')} {Z^2}_{\hspace{-.6cm}i\; i',22 \; 22'} \ , \\
\ D_{23} =& \prod\limits^{17}_{\stackrel{i=1}{i \neq 8-10,12}}
\stackrel{\scriptstyle (19)-(22')}{Z^2}_{\hspace{-.6cm}i\; i',23 \; 23'} \ ,
\ D_{24} = \prod\limits^{18}_{\stackrel{i=1}{i \neq 11,12}}
\stackrel{\scriptstyle (19)-(23')}{Z^2}_{\hspace{-.6cm}i\; i',24 \; 24'}.
\end{align*}
Let us denote the regenerations of $\tilde{C}_i$  as $C_i$, $1 \leq i \leq 10$.
These are the factorizations of the suitable $D_t$, as in (\ref{formulac}).

\subsection{Properties of the factorization $\Delta_{48}^2$}\label{prop}
The braid monodromy factorization $\Dl_{48}^2$
is the product $\prodl_{i=1}^{10}C_i \varphi_i$.
We now verify that there are no missing braids.

\begin{proposition}
The braid monodromy factorization of $S$ is $\Delta_{48}^2$.
\end{proposition}

\begin{proof}
By Proposition  V.2.2 in \cite{16},  deg$\Dl_{48}^2=48 \cdot 47=2256$.

Now we check the degree of the factorization $\prodl_{i=1}^{10}C_i \varphi_i$.
The monodromies $\set{\varphi_i}_{i=1, 3, 6, 8}$ each consist of six cusps and
nine branch points, see Section \ref{mono3}.
We combine their degrees to get $4\cdot(6 \cdot 3 + 1 \cdot 9)=4\cdot 27=108$.

In each $\varphi_i, i=2, 4, 5, 7, 9, 10$ (Section \ref{mono6}), $\deg F(F)^{*} = 48$.
The factors outside  $F(F)^{*}$ are $20$ degree two factors,
$12$ degree three factors, and $2$ degree one factors.
Combining the degrees of these,
we get $20 \cdot 2 + 12 \cdot 3 + 2 \cdot 1 = 78$.
Thus $\deg\varphi_i = 78 + 48 = 126$ for $i= 2, 4, 5, 7, 9, 10$,
and therefore $\deg(\prodl_{i= 2, 4, 5, 7, 9, 10}\varphi_i) = 126 \cdot 6 = 756$.

The products $\set{C_i}_{i=1}^{10}$ consist of $696$ nodes, see Section \ref{regci}.
We combine degrees of their factors to get
$\deg\prodl_{i=1}^{10}{C_i}=696\cdot2=1392$.

Finally $\deg\prodl^{10}_{i=1}C_i \varphi_i = 108+756+1392 = 2256$.
Therefore $\Delta_{48}^2$ is the desired braid monodromy factorization.
\end{proof}

\medskip

We now study invariance properties of $\Delta^2_{48}$.
Invariance properties are results
in which we prove that the braid monodromy factorization $\Delta^2_{48}$
is invariant under certain elements of $B_{48}$.
Establishing invariance properties is essential
in order to simplify the computations
which follow from the van Kampen Theorem (see \cite{pil});
at the end of this section
we describe the effect of these properties on the application of the van Kampen Theorem.

The following definitions are necessary in order to prove the
Invariance Properties.

\begin{definition}\label{df:154}
Let $g_1 \cdots g_k = h_1  \cdots  h_k$ be two factorized expressions
of the same element in a group $G$.
We say that $g_1 \cdots g_k $ is obtained from $h_1 \cdots h_k$
by a Hurwitz move if $\exists \ 1 \leq p \leq k- 1$,
such that $g_i = h_i \ (i \neq p,p + 1) \  g_p = h_p h_{p+1}h^{-1}_p$
and $g_{p+1} = h_p$ or $g_p = h_{p+1}$ and
$g_{p+1} = h^{-1}_{p+1} h_p h_{p+1}$.
In general,  $g_1 \cdots  {g_k}  \simeq  h_1 \cdots h_k$ (Hurwitz equivalent) if
$g_1 \cdots  {g_k}$ is obtained from $h_1 \cdots  h_k$
by a finite number of Hurwitz moves.
\end{definition}

\begin{definition}\label{df:156}
Let $g_1 \cdots g_k$  be a factorized expression in $G$ and $h \in G$.
We say that $g_1 \cdots g_k$ is invariant under $h$
if it is Hurwitz equivalent to $(g_1)_h \cdots (g_k)_h$,
where $(g_i)_h = h^{-1} g_ih$.
\end{definition}

Invariance properties are important in view of Lemma VI.4.2 in \cite{16}:
if a braid monodromy factorization $\Delta^2_{n}$ is invariant under $h$
then the equivalent factorization $(\Delta^2_{n})_h$
is also a braid monodromy factorization.

We now quote several invariance rules
(the conclusions of Theorems \ref{1rule}, \ref{2rule}, \ref{3rule}),
Chakiri's Lemma (see \cite{19}), and make some remarks on invariance properties.

\underline{Invariance Rules}
\begin{enumerate}\label{invrul}
\item
A braid $Z^{1}_{ij}$ is invariant under $(Z_{ii'} Z_{jj'})^q$ for any $q \in \Z$.
\item
A braid $Z^{2}_{i,jj'}$ (resp. $Z^{2}_{ii',jj'}$) is invariant under $Z_{jj'}^q$
(resp. $Z_{ii'}^p Z^q_{jj'}$) for any $p, q \in \Z$.
\item
$Z^{3}_{i,jj'}$ is invariant under $Z_{jj'}^q$ for any $q \in \Z$.
\end{enumerate}

\begin{lemma}\label{lm:159}{\bf{(Chakiri)}.}
Let $g = g_1 \cdots g_k$ be a factorized expression in a group $G$.
Then $g_1 \cdots g_k$ is invariant under $g^m$ for any $m \in \Z$.
\end{lemma}

\underline{Invariance Remarks}
\begin{enumerate}
\item
To prove invariance of $g_1 \cdots g_k$ under $h$
it is enough to prove that $g_1 \cdots g_t$
and $g_{t+1} \cdots g_k$ are invariant under $h$.
Thus we can divide a factorization into
subfactorizations and prove invariance on each part separately.
\item
An element $g_1$ is invariant under $h$ if and only if $g_1$ commutes with $h$.
\item
If a product of elements that commutes with $h$ is invariant under $h$,
the corresponding factorizations are equal.
\item
If two paths $\sigma_1$ and $ \sigma_2$ do not intersect,
the corresponding halftwists $H(\sigma_1)$ and $H(\sigma_2)$ commute.
\item
If $g$ is invariant under $h_1$ and $h_2$ then $g$ is invariant under $h_1 h_2$.
\end{enumerate}

We finally prove the Invariance Properties of $\Delta^2_{48}$.
Denote $Z_{ij} = H(z_{ij})$.

\begin{remark}\label{rem-ij}
The halftwists $Z_{ii'}$ and $Z_{jj'}$ commute for all $i, j$,
since the path from $i$ to $i'$ does not intersect the path from $j$ to $j'$.
\end{remark}

\begin{lemma}\label{inva1357}
Each monodromy among $\varphi_1, \varphi_3, \varphi_6, \varphi_8$
is invariant under  $\prodl^{24}_{j=1} Z^{m_j}_{jj'}$  for $m_j \in \Z$.
\end{lemma}

\begin{proof}
The braids in $\varphi_1, \varphi_3, \varphi_6, \varphi_8$, arise from cusps and branch points
(see Theorem \ref{second2}).
By Invariance Rule (3),
each braid of the form $Z^{3}_{ii',j}$ is invariant under $Z_{ii'}^q$.
Each of the braids of the form $Z^{1}_{ij}$
is invariant under $(Z_{ii'} Z_{jj'})^q$, and in particular,
if the braid is of the form  $Z^{1}_{ii'}$, then
it is invariant under $Z_{ii'}^q$.

Moreover, by Remark \ref{rem-ij},
each $Z_{kk'}$ commutes with the above braids when $k \neq i, j$.

Therefore, each monodromy $\varphi_i$, $i=1, 3, 6, 8$,
is invariant under $\prodl^{24}_{j=1} Z^{m_j}_{jj'}$ for $m_j \in \Z$.
\end{proof}

\begin{lemma}\label{invaci}
Each factorization $\{C_i\}$, $i=1, \dots, 10$, is invariant under
$\prodl^{24}_{j=1} Z^{m_j}_{jj'}$  for $m_j \in \Z$.
\end{lemma}

\begin{proof}
We apply Invariance Rule (2) and Invariance Remark (4) to each $C_i$
to get the desired invariance.
\end{proof}

We now prove invariance for the monodromy $\varphi_2$.
A similar proof will apply for each one of the monodromies
$\varphi_4, \varphi_5, \varphi_7, \varphi_9, \varphi_{10}$.

\begin{lemma}\label{invaphi2}
$\varphi_2$ is invariant under $(Z_{4\; 4'} Z_{5' \; 5'})^p (Z_{13 \;13'}Z_{20 \;20'})^q
(Z_{1\; 1'} Z_{2\; 2'})^r$ for all $p,q,r\in \Z$.
\end{lemma}

\begin{proof}

\underline{Case 1}: $p = q = r$.\\
As proved in \cite[Lemma 12]{19}, the monodromy $\varphi_2$ can be written as
$\prodl^{}_{j=1,2,4,5,13,20} Z_{jj'}^{-1} \Delta_{12}^2$.
By Lemma \ref{lm:159}, $\varphi_2$
is invariant under  $(\prodl^{}_{j=1,2,4,5,13,20} Z_{jj'}^{-1} \Delta_{12}^2)^{-p}$.
Since $\Delta_{12}^2$ is a central element,
$\varphi_2$ is invariant under
$(Z_{4\; 4'} Z_{5 \; 5'})^p (Z_{13 \; 13'}Z_{20 \; 20'})^p (Z_{1\; 1'} Z_{2\; 2'})^p$.

\underline{Case 2}: $p = 0$.\\
Let $\epsilon = (Z_{13 \;13'} Z_{20 \;20'})^q (Z_{1\; 1'}Z_{2\; 2'})^r$.
\begin{enumerate}
\item
\underline{Step 1:}  Factors outside of $F \cdot F^{Z_{20 \; 20'}^{-1}Z_{13 \; 13'}^{-1}}$. \\
${Z_{4 \; 4'}}^{Z^2_{4',13 \; 13'}Z^2_{4',5 \;5'}Z^2_{2 \; 2',4}Z^2_{1 \; 1',4}}$ and
${Z_{5 \; 5'}}^{Z^2_{2 \; 2',5}Z^2_{5',20 \;20'}Z^2_{5',13 \; 13'}Z^2_{2 \; 2',4}Z^2_{1 \; 1',4}}$
commute with $\epsilon$.\\
$Z^3_{1 \;1',4}, {Z^3_{4',13 \;13'}}^{Z^{2}_{4' ,5 \; 5'}}, {Z^3_{2 \;2',5}}^{Z^2_{2 \; 2',4}Z^2_{1 \; 1',4}}$
and ${Z^3_{5',20 \;20'}}^{Z^{2}_{5' ,13 \; 13'}Z^2_{2 \; 2',4}Z^2_{1 \; 1',4}}$
are invariant under $\epsilon$ by Invariance Rule (3). \\
The degree 2 factors are of the form $Z^2_{\alpha \alpha ', \beta}$ where $\beta = 4,4',5,5'$.
By Invariance Rule (2), they are invariant under $Z_{\alpha \alpha '}$,
and since the other halftwists in $\epsilon$ commute with $Z_{\alpha \beta}$ and
$Z_{\alpha ' \beta }$,
we get that $Z^2_{\alpha \alpha ', \beta }$ is invariant under $\epsilon$.\\
Moreover all conjugations which appear in $\varphi_2$ (i.e. ${( \ \ )}^*$)
are invariant under $\epsilon$
by Invariance Rule (2) and by Invariance Remark (4).
\item
\underline{Step 2}:  Factors in  $F \cdot F^{Z_{20 \; 20'}^{-1}Z_{13 \; 13'}^{-1}}$.\\
Let $\rho = Z_{13 \; 13'} Z_{20 \; 20'}$.
In order to prove that  $F \cdot F^{\rho^{-1}}$ is invariant under $\epsilon$,
we consider the following subcases:
\begin{itemize}
\item
\underline{Subcase 2.1}:  $q = 0$  and $\epsilon = (Z_{1\; 1'}Z_{2\; 2'})^r$.\\
$Z^3_{2 \; 2',13}$ and ${Z^3_{2 \; 2',20}}^{Z^2_{2 \; 2',13}}$
are invariant under $Z_{2\; 2'}$ (Invariance Rule (3)) and commute with $Z_{1\; 1'}$
(Invariance Remark (4)).
Thus $Z^3_{2 \; 2',13}{Z^3_{2 \; 2',20}}^{Z^2_{2 \; 2',13}}$ is
invariant under $\epsilon$ (Invariance Remarks (1) and (5)).
$Z^2_{13' \; 20}$ and $\cdot(Z^2_{13 \; 20})^{Z^2_{2 \; 2',13}}$ commute with
$Z_{1\; 1'}$ and $Z_{2\; 2'}$ and thus with $\epsilon$.
$\tilde{Z}_{1 \; 2'} \cdot \tilde{Z}_{1' \; 2}$ is invariant under
$(Z_{1\; 1'}  Z_{2\; 2'})^r$ by Invariance Rule (1).
For $F^{\rho^{-1}}$ we have the same conclusion,
since $Z^{-1}_{20 \; 20'} Z^{-1}_{13 \;13'}$ commutes with $\epsilon = (Z_{1\; 1'}Z_{2\; 2'})^r$.
\item
\underline{Subcase 2.2}:  $r=0 \ , \ q = 1$  and  $\epsilon = Z_{13 \;13'}Z_{20 \;20'}$.\\
Note that in this case $\epsilon = \rho$.
To prove that  $F \cdot F^{\rho^{-1}}$  is invariant under $\rho$,
we must show that $F \cdot F^{\rho^{-1}}$ is Hurwitz equivalent to
$F^{\rho} \cdot F$.
Since $AB$ is Hurwitz equivalent to $BA^B$,
it is enough to prove that   $F \cdot F^{\rho^{-1}}$ is Hurwitz equivalent to
$F \cdot {(F^{\rho})}^F$.
Thus it is enough to prove that
$F^{\rho^{-1}}$ is Hurwitz equivalent to ${(F^{\rho})}^F$
or that ${(F^{\rho^{-1}})}^{F^{-1}}$ is Hurwitz equivalent to $F^{\rho}$.

By Theorem \ref{mono2},
$F = \Delta^2_8 \rho^{-2} Z^{-2}_{1\; 1'} Z^{-2}_{2\; 2'} (F^{-1})^{\rho^{-1}}$,
thus $F^{-1} = F^{\rho^{-1}} (Z_{1\; 1'} Z_{2\; 2'})^2 \rho^2 \Delta^2_8$.
Now we have
\begin{eqnarray*}
(F^{\rho^{-1}})^{F^{-1}}  =  (F^{\rho^{-1}})^{F^{\rho^{-1}} (Z_{1\; 1'}Z_{2\; 2'})^2 \rho^2 \Delta^2_8}
\left . \begin{array}{c} = \\
[-.5cm]{\scriptstyle as \ factorized} \\ [-.5cm]{\scriptstyle expression}
\end{array}  \right .
(F^{\rho^{-1}})^{F^{\rho^{-1}} (Z_{1\; 1'}Z_{2\; 2'})^2 \rho^2}\\
\left . \begin{array}{c}
\simeq \\ [-.5cm]{\scriptstyle Chakiri}  \end{array}  \right .
(F^{\rho^{-1}})^{(Z_{1\; 1'}Z_{2\; 2'})^2 \rho^2} =  {(F^{(Z_{1\; 1'}Z_{2\; 2'})^2})}^\rho \; \left .
\begin{array}{c}
\simeq \\ [-.5cm]{\scriptstyle \ Subcase \ 2.1}
\end{array}  \right . F^{\rho}.
\end{eqnarray*}
\item
\underline{Subcase 2.3}:  $q=2q' \ ; \ \epsilon = (Z_{13 \;13'}Z_{20 \;20'})^{2q'}
(Z_{1\; 1'}Z_{2\; 2'})^r$. \\
$F \cdot F^{\rho^{-1}}$ can be written as
$\Delta^2_8(Z_{13 \;13'} Z_{20 \;20'})^{-2}(Z_{1\; 1'}Z_{2\; 2'})^{-2}$.
By Chakiri's Lemma,
$F \cdot F^{\rho^{-1}}$ is invariant under
${(\Delta^2_8(Z_{13 \;13'} Z_{20 \;20'})^{-2}(Z_{1\; 1'}Z_{2\; 2'})^{-2})}^{-q'}$ and thus under
$(Z_{13 \;13'}Z_{20 \;20'})^{2q'}(Z_{1\; 1'}Z_{2\; 2'})^{2q'}$.
By Subcase 2.1,  $F \cdot F^{\rho^{-1}}$
is invariant under $(Z_{1\; 1'}Z_{2\; 2'})^{r-2q'}$.
By Invariance Remark (5), $F \cdot F^{\rho^{-1}}$
is invariant under $$(Z_{13 \;13'}Z_{20 \;20'})^{2q'} (Z_{1\; 1'}Z_{2\; 2'})^{2q'}
(Z_{1\; 1'}Z_{2\; 2'})^{r-2q'} = \epsilon.$$
\item
\underline{Subcase 2.4}:  $q=2q'+1 \ ; \  \epsilon = (Z_{13 \;13'}Z_{20 \;20'})^{2q'+1}
(Z_{1\; 1'}Z_{2\; 2'})^r$. \\
It is easy to verify this case by using the prior cases 2.2, 2.3 and Invariance Remark (5).
\end{itemize}
\end{enumerate}

\underline{Case 3}:  $p,q,r$ arbitrary; $\epsilon = (Z_{4\; 4'} Z_{5 \; 5'})^p (Z_{13 \;13'}
Z_{20 \;20'})^q (Z_{1\; 1'} Z_{2\; 2'})^r$. \\
By case 1, $\varphi_2$ is invariant under $(Z_{4\; 4'} Z_{5 \; 5'})^p (Z_{13 \;13'}
Z_{20 \;20'})^p (Z_{1\; 1'} Z_{2\; 2'})^p$, and by case 2 under
$(Z_{13 \;13'}Z_{20 \;20'})^{q-p} (Z_{1\; 1'} Z_{2\; 2'})^{r-p}$.
By Invariance Remark (5),
$\varphi_2$ is invariant under $\epsilon$.
\end{proof}

\begin{corollary}\label{inva6pt}
Each one of the monodromies
$\varphi_4, \varphi_5, \varphi_7, \varphi_9, \varphi_{10}$
is invariant under \\ $(Z_{ii'} Z_{jj'})^p (Z_{kk'}Z_{ll'})^q
(Z_{mm'} Z_{nn'})^r$ for all $p,q,r \in \Z$,
where $i$ and $j$ are the diagonal lines
around the relevant $6$-point,
and $k$ and $l$ (resp. $m$ and $n$) are the vertical
(resp. horizontal) ones (see Figure \ref{top}).
\end{corollary}

Since the proof of the invariance for the $6$-points
relies on the invariance of each local braid monodromy
under such expressions as in Corollary \ref{inva6pt},
we have to pay attention that each diagonal line
(except the lines 14, 17, 19, 24 in Figure \ref{top})
at a certain $6$-point is also a diagonal line at another $6$-point.
For example, the monodromy $\varphi_2$ is invariant under
$(Z_{4\; 4'} Z_{5 \; 5'})^p (Z_{13 \;13'}Z_{20 \;20'})^q (Z_{1\; 1'} Z_{2\; 2'})^r$ for all $p,q,r\in \Z$
(see Lemma \ref{invaphi2}).
In addition, the monodromy $\varphi_4$ is invariant under
$(Z_{4\; 4'} Z_{8\; 8'})^{p'} (Z_{3\; 3'}Z_{7\; 7'})^{q'}
(Z_{15 \;15'} Z_{21 \;21'})^{r'}$ for all $p',q',r' \in \Z$ (see Corollary \ref{inva6pt}).
This implies that $p=p'$.
In the same way,
$\varphi_5$ is invariant under  $(Z_{5\; 5'} Z_{9 \; 9'})^{p''} (Z_{6 \; 6'}Z_{11 \; 11'})^{q''}
(Z_{16\; 16'} Z_{22\; 22'})^{r''}$ for all $p'',q'',r''\in \Z$.
This implies that $p=p''$.
In the notation of Lemmas \ref{inva1357} and \ref{invaci},
$p=m_{4}=m_{5}, p'=m_{4}=m_8, p''=m_5=m_{9}$, and
therefore $m_{4} =  m_{5} =  m_{8} =  m_{9}$.
Now, the lines 14 and 17 (resp. 19 and 24) are diagonal lines at only one $6$-point 9 (resp. 10).
Therefore in the above notation, $m_{14} = m_{17}$ (resp. $m_{19} = m_{24}$).

We now consider the vertical lines.
By a similar argument as above, we can conclude that
$m_{13} = m_{18} = m_{20} = m_{23}$
from the invariance of $\varphi_2, \varphi_7, \varphi_9$ and $\varphi_{10}$.
Since the lines $6$ and $11$ (resp. $3$ and $7$)
are vertical lines at only one $6$-point 5 (resp. 4),
we get $m_{6} = m_{11}$ (resp. $m_{3} = m_{7}$).

For the horizontal lines,
we can immediately conclude that $m_{1} = m_{2}, m_{10} = m_{12}$, and
$m_{15} = m_{16} = m_{21} = m_{22}$.

\begin{corollary}\label{finalinva}
The braid monodromy factorization $\Delta_{48}^2$ is invariant under
\begin{eqnarray*}
\rho = &{(Z_{1\; 1'}Z_{2\; 2'})}^{m_1} {(Z_{6\; 6'}Z_{11\; 11'})}^{m_6} {(Z_{10\; 10'}Z_{12\; 12'})}^{m_{10}}
{(Z_{3\; 3'}Z_{7 \;7'})}^{m_3} {(Z_{14 \;14'}Z_{17 \;17'})}^{m_{14}} {(Z_{19 \;19'}Z_{24 \;24'})}^{m_{19}} \\
& {(Z_{4\; 4'}Z_{5\; 5'}Z_{8 \; 8'}Z_{9 \; 9'})}^{m_{4}}
{(Z_{13 \;13'}Z_{18 \; 18'}Z_{20 \; 20'}Z_{23 \; 23'})}^{m_{13}}
{(Z_{15 \;15'}Z_{16 \; 16'}Z_{21 \; 21'}Z_{22 \; 22'})}^{m_{15}}.
\end{eqnarray*}
\end{corollary}

\bigskip

\subsection{Consequences from the Invariance Theorems}\label{imp}
We are interested in the fundamental group $\pi_1(\C\P^2 - S)$
and the fundamental group of the Galois cover of the surface
(with respect to the generic projection onto $\C\P^2$).
These groups are computed in \cite{pil}.

Here we extract qualitative information
concerning the importance of the above computations and
how they are connected to the groups.

\begin{figure}[ht]
\begin{minipage}{\textwidth}
\begin{center}
\epsfbox{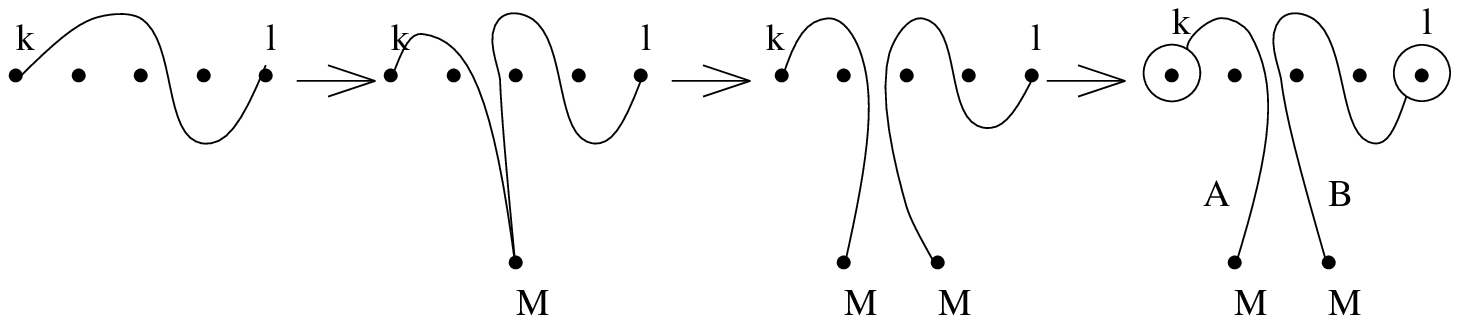}
\end{center}
\end{minipage}
\caption{}\label{AB}
\end{figure}

\begin{figure}[ht!]
\epsfxsize=7cm 
\epsfysize=6cm 
\begin{minipage}{\textwidth}
\begin{center}
\epsfbox{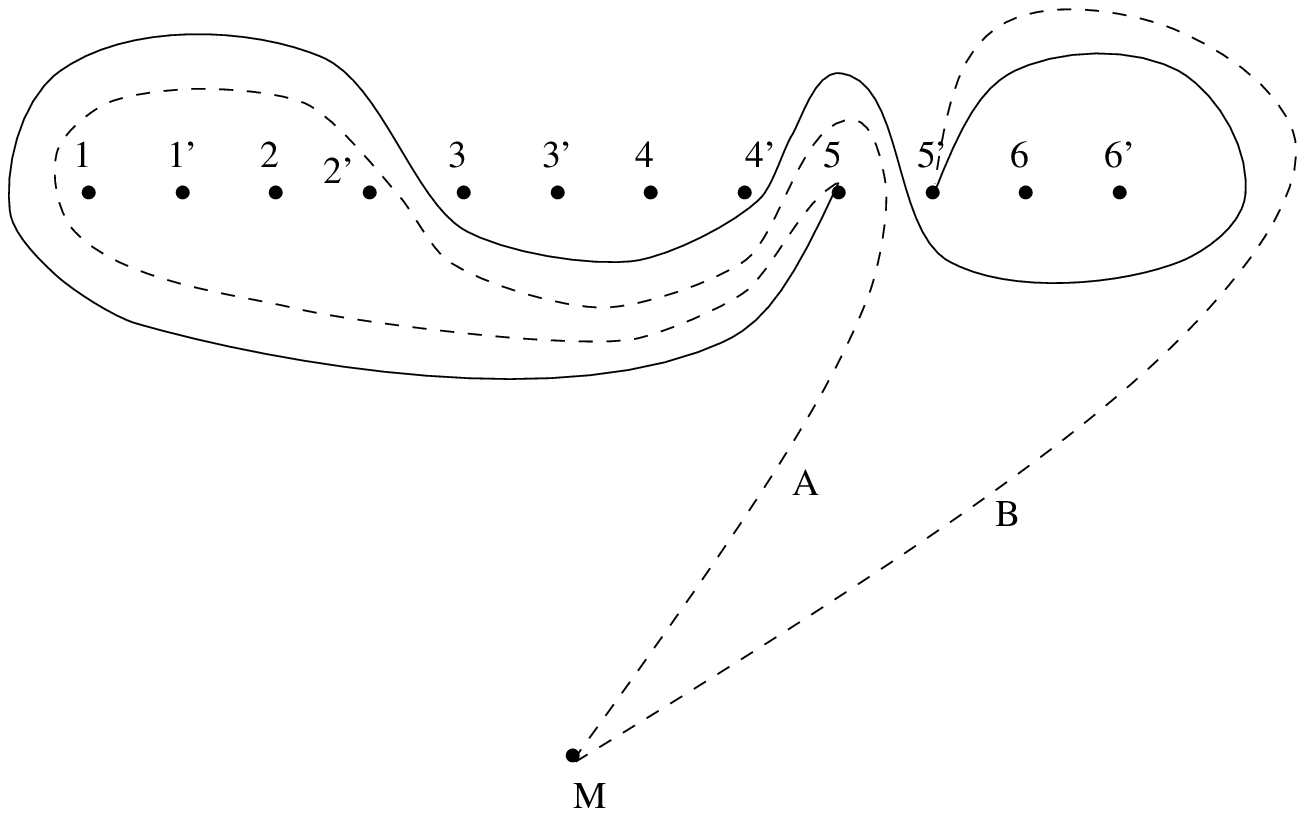}
\end{center}
\end{minipage}
\caption{}\label{vk1}
\end{figure}

The van Kampen Theorem \cite{20}
induces a finite presentation of the fundamental group of complements of curves
by means of generators and relations.
This is done by applying the theorem on the resulting factorization $\Delta^2_{48}$.
We take any path from $k$ to $\ell$, cut it in $M$, then towards $k$ along the path,
around $k$ and coming back the same way.
Consider $A$ as an element of the fundamental group.
Do the same to $\ell$ to obtain $B$.
See Figure \ref{AB}.
Now, according to the van Kampen Theorem we obtain a relation which involves  $A$ and  $B$.
The following example is illustrative.
Consider Figure \ref{vk1}.
We construct $B$ by proceeding from $M$ towards $5'$ above $6'$ and 6 encircling $5'$ counterclockwise
and proceeding back above 6 and $6'$.
Therefore $B = {6'} 6 {5'} 6^{-1} {6'}^{-1}$.
In a similar way,
$A = 5 {2'} 2 {1'} 1 5 1^{-1} {1'}^{-1} 2^{-1} {2'}^{-1} 5^{-1}$.
If this path is obtained from a branch point (by Moishezon-Teicher monodromy),
then by the van Kampen Theorem, we have the relation
${6'} 6 {5'} 6^{-1} {6'}^{-1} = 5 {2'} 2 {1'} 1 5 1^{-1} {1'}^{-1} 2^{-1} {2'}^{-1} 5^{-1}$.

\begin{figure}[ht!]
\epsfxsize=8cm 
\epsfysize=6cm 
\begin{minipage}{\textwidth}
\begin{center}
\epsfbox{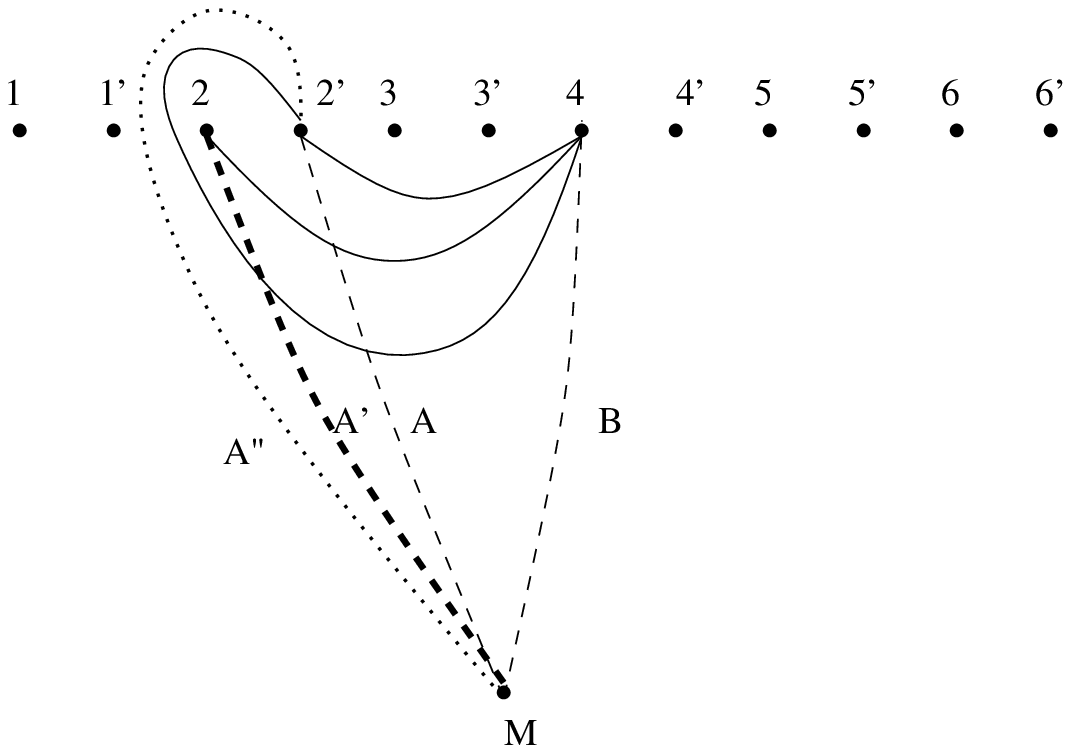}
\end{center}
\end{minipage}
\caption{}\label{vk2}
\end{figure}

\begin{figure}[ht!]
\epsfxsize=7cm 
\epsfysize=6cm 
\begin{minipage}{\textwidth}
\begin{center}
\epsfbox{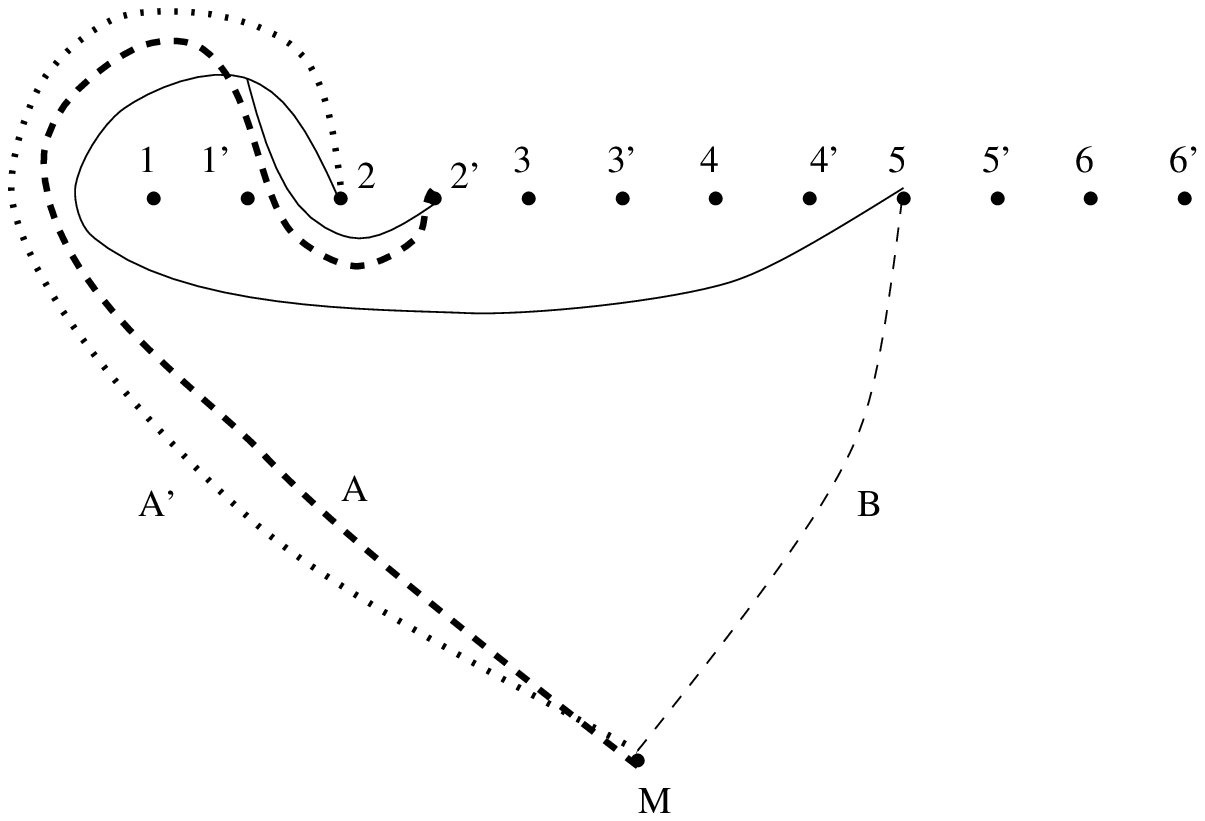}
\end{center}
\end{minipage}
\caption{}\label{vk3}
\end{figure}

In a similar way, consider paths which are related to cusps (Figure \ref{vk2}) and nodes
(Figure \ref{vk3}).
From Figure \ref{vk2} we get three relations $B A B = A B A$,
$B A' B = A' B A'$ and $B A'' B = A'' B A''$.
And from Figure \ref{vk3}, we get $A B = B A$ and $A' B = B A'$.

We now explain the importance of the invariance of $\Delta_{48}^2$ (Corollary \ref{finalinva}).
Let $z_{ij}$ be a path connecting $i$ or
${i'}$ with $j$ or ${j'}$ and $Z_{ij}$ its corresponding halftwist.
We can conjugate  $Z_{ij}$ by $Z_{ii'}^{\pm 1}, Z_{jj'}^{\pm 1}$.
These conjugations are the actions of $Z_{ii'}^{\pm 1}$ (resp. $Z_{jj'}^{\pm 1}$)
on the ``head'' (resp. ``tail'') of $z_{ij}$ within a small circle around $i$ and ${i'}$
(resp. $j$ and $j'$).
The ``body'' of $z_{ij}$ does not change under such conjugations
and in particular not under $Z_{ii'}^{\pm m}$ or $Z_{jj'}^{\pm n}$  for $m, n \in \Z$.
Therefore, Corollary \ref{finalinva} enables us to use Theorem 1.6 from \cite{MoTe10}:

\begin{theorem}
If a sub-factorization $\prodl^r_{i=s} Z_i$ is invariant
under any element $h$,  and $\prodl^r_{i=s} Z_i$ induces a relation
$A_{i_1} \cdot ... \cdot A_{i_t}$ on  $\pi_1(\C\P^2 - S)$
via the van Kampen Theorem, then $(A_{i_1})_h  \cdot  ... \cdot (A_{i_t})_h$ is also a relation.
\end{theorem}

That means that we can expand our list of relations,
which assists us in the computations in \cite{pil}:

\begin{corollary}
Consider $\rho =  \prodl^{24}_{j=1} \rho_{m_j}= \prodl^{24}_{j=1} Z^{m_j}_{jj'}$
from Corollary \ref{finalinva}.
If $R$ is any relation in $\pi_1(\C\P^2 - S)$,
then $R_\rho$ is also a relation in  $\pi_1(\C\P^2 - S)$,
where $R_\rho$ is the relation induced from
$R$ by replacing $A_j$ with $(A_j)_{\rho_{m_j}}$.
\end{corollary}

\section{Acknowledgements}

This research was initiated while the first author was staying at the Mathematics
Institute, Erlangen - N\"urnberg university, Germany.
She wishes to thank the Institute for their hospitality
and especially acknowledge her hosts Wolf Barth and Herbert Lange.
She wishes to thank also the Einstein Institute for Mathematics (Jerusalem)
for her present stay, and to her hosts Hershel Farkas and Ruth Lawrence-Neumark.

Thanks are given to Michael Friedman for several fruitful discussions.

\end{document}